\documentclass[11pt,a4paper, reqno]{amsart}
\usepackage{amsmath,amssymb}
\usepackage{bigints}
\usepackage{indentfirst}
\usepackage{fancybox}
\usepackage{ascmac}
\usepackage{latexsym}
\usepackage{bm}
\usepackage{amsfonts}
\usepackage[all]{xy}
\usepackage{xcolor}
\usepackage{tikz}
\usetikzlibrary{cd}
\usetikzlibrary{intersections, calc, arrows}
\usepackage{here}
\usepackage{amsthm}
\usepackage{comment}

\usepackage{mathtools}
\mathtoolsset{showonlyrefs,showmanualtags} 

\makeatletter
\@addtoreset{equation}{section}
\makeatother

\makeatletter 
\newcommand\fleqnoff{\@fleqnfalse\@mathmargin\@centering} 
\newcommand\fleqnon[1][\leftmargini]{\@fleqntrue\@mathmargin=#1\relax 
\@ifundefined{mathindent}{\let\mathindent\@mathmargin}{}} 
\makeatother

\title[]{On Triple Quadratic Residue Symbols in Real Quadratic Fields}
\author{Atsuki Kuramoto}
\address{Graduate School for Mathematics, Kyushu University, Motooka 744, Nishi-ku Fukuoka 819-0395, Japan}
\email{kuramoto.atsuki.2358@gmail.com}
\date{\today}
\keywords{Pro-2 Galois group with restricted ramification, Mod 2 arithmetic Milnor invariant, R{\'e}dei extension, Triple quadratic residue symbol}
\subjclass[2020]{11R32, 57M05}

\theoremstyle{definition}

\newtheorem{thm}{Theorem}[section]
\newtheorem*{thm*}{Theorem}
\newtheorem{dfn}[thm]{Definition}
\newtheorem{lem}[thm]{Lemma}
\newtheorem{prp}[thm]{Proposition}
\newtheorem{cor}[thm]{Corollary}
\newtheorem{ex}[thm]{Example}
\newtheorem{rem}[thm]{Remark}

\newcommand{\Thm}[3]{ 
\begin{thm}[#1]\label{#2}
 \textit{#3}
\end{thm}}

\newcommand{\Def}[3]{
 \begin{dfn}[#1]\label{#2}
  #3
 \end{dfn}}

\newcommand{\Lem}[3]{
\begin{lem}[#1]\label{#2}
  \textit{#3}
 \end{lem}}

\newcommand{\Prp}[3]{
 \begin{prp}[#1]\label{#2}
  \textit{#3}
 \end{prp}}

\newcommand{\Cor}[3]{
 \begin{cor}[#1]\label{#2}
  \textit{#3}
 \end{cor}
}

\newcommand{\Rem}[3]{
\begin{rem}[#1]\label{#2}
 #3
\end{rem}}

\newcommand{\Ex}[3]{
\begin{ex}[#1]\label{#2}
 #3
\end{ex}}

\begin{document}

\begin{abstract}
We introduce triple quadratic residue symbols $[\frak{p}_1, \frak{p}_2, \frak{p}_3]$ for certain finite primes $\frak{p}_i$'s of a real quadratic field $k$ with trivial narrow class group.
For this, we determine a presentation of the Galois group of the maximal pro-2 Galois extension over $k$ unramified outside $\frak{p}_1, \frak{p}_2, \frak{p}_3$ and infinite primes, from which we derive mod 2 arithmetic triple Milnor invariants $\mu_2(123)$ yielding the triple symbol $[\frak{p}_1, \frak{p}_2, \frak{p}_3] = (-1)^{\mu_2(123)}$.
Our symbols $[\frak{p}_1, \frak{p}_2, \frak{p}_3]$ describes the decomposition law of $\frak{p}_3$ in a certain dihedral extension $K$ over $k$ of degree 8, determined by $\frak{p}_1, \frak{p}_2$.
The field $K$ and our symbols $[\frak{p}_1, \frak{p}_2, \frak{p}_3]$ are generalizations over real quadratic fields of R\'{e}dei's dihedral extension of $\mathbb{Q}$ and R\'{e}dei's triple symbol of rational primes.
We give examples of R\'{e}dei type extensions $K$ over real quadratic fields.
We also give a cohomological interpretation of our symbols in terms of Massey products.
\end{abstract}

\maketitle

\section*{Introduction}
The purpose of this paper is to introduce the triple quadratic residue symbol $[\frak{p}_1, \frak{p}_2, \frak{p}_3]$ for certain primes $\frak{p}_i$'s of real quadratic fields and investigate its properties.
Our symbols is a triplication of quadratic residue symbol in real quadratic fields, and is also a generalization of R\'{e}dei's triple symbol over the rational number field $\mathbb{Q}$ to real quadratic fields.

The study of triple symbols in arithmetic goes back to the work of R\'{e}dei in 1939 (\cite{Redei1939}), which intended to generalize the Legendre symbol and Gauss' genus theory for quadratic fields.
For certain prime numbers $p_1, p_2$ and $p_3$, R\'{e}dei's triple symbol $[p_1, p_2, p_3] \in \{ \pm 1\}$ describes the decomposition law of $p_3$ in the dihedral extension of degree 8 over $\mathbb{Q}$, determined by $p_1$ and $p_2$.
Some variants of the R\'{e}dei symbol have also been studied in \cite{frohlich1960prime}, \cite{Furuta1980}, \cite{Suzuki} among others.
In the late 1990s, Morishita interpreted the R\'{e}dei symbol as a mod 2 arithmetic analogue of the triple linking number (Milnor invariant) of a link $L$ in the 3-sphere $S^3$ (\cite{Milnor1957}, \cite{Turaev1979}) from the viewpoint of arithmetic topology (\cite{Morishita2024}), especially from the analogy between Koch's theorem on a group presentation of the Galois group $\pi_1^{\footnotesize \mbox{pro-2}}({\rm Spec}(\mathbb{Z}) \setminus S)$ of the maximal pro-2 extension over $\mathbb{Q}$ unramified outside a finite set $S$ of primes and Milnor's theorem on a group presentation of a link group $\pi_1(S^3 \setminus L)$, and introduced the multiple quadratic residue symbols $[p_1, \dots , p_r]$ for certain prime numbers $p_i$'s as mod 2 arithmetic analogues of Milnor invariants of a link (\cite{Morishita2000}, \cite{Morishita2002}, \cite{Morishita2004}).
He also gave a description of $[p_1, \dots , p_r]$ in terms of the Massey product of the \'{e}tale cohomology of ${\rm Spec}(\mathbb{Z}) \setminus \{(p_1), \dots, (p_r)\}$ (\cite{Morishita2004}).

Now it remains an important problem to generalize the multiple quadratic residue symbols in $\mathbb{Q}$ to the multiple power residue symbols for primes in a number field. We may note that Turaev (\cite{Turaev1979}) extended Milnor's theorem on a link group to an integral homology $3$-sphere $M$ and defined Milnor's invariants of a link in $M$.
Since the ideal class group is regarded as an analogue of the $1$st homology group of a $3$-manifold  in arithmetic topology, it is natural to expect that we may extend Koch's presentation of the pro-$2$ Galois group over $\mathbb{Q}$ to the Galois group $G_{k, S}(l)$ of the maximal pro-$l$ extension over a number field $k$ unramified outside a finite set $S$ of primes, when $k$ has the (narrow) class number one and $k$ contains the $l$-th roots of unity ($l$ being a prime number), and that we may derive multiple $l$-th power residue symbols from $G_{k, S}(l)$.
However, there are a cohomological obstruction to determine a group presentation of $G_{k, S}(l)$ and the difficulty arising from the unit group of $k$, in order to define arithmetic Milnor invariants.
So far there is the only case $k = \mathbb{Q}(\sqrt{-3})$ and $l=3$, for which triple cubic residue symbols are defined by using the Galois group $G_{\mathbb{Q}(\sqrt{-3}), S}(3)$ (\cite{AMM}).
For other approaches to triple symbols, we refer to the recent works \cite{Efrat2024} and \cite{2407.02063}.

In this paper, we resolve the difficulty arising from the cohomological obstruction and the unit group for a real quadratic field $k$ and determine a presentation of the Galois group $G_S(2)$ of the maximal pro-2 extension of $k$ unramified outside a finite set of finite primes $S = \{\frak{p}_1, \dots \frak{p}_s\}$ and the infinite primes, when $k$ has the trivial narrow class group.
From this we are able to derive the triple quadratic residue symbol $[\frak{p}_1, \frak{p}_2, \frak{p}_3] = \pm1$ for certain finite primes of $k$.
To be precise, fix a prime number $p \equiv 1$ mod 4 and let $k$ be a real quadratic field $\mathbb{Q}(\sqrt{p})$ whose narrow ideal class group is trivial, and let $\varepsilon$ be the fundamental unit of $k$.
Let $S = \{ \frak{p}_1, \dots, \frak{p}_s \}$ be a finite set of finite primes of $k$ and let $\infty_1, \infty_2$ denotes the infinite primes of $k$, where ${\rm N}\frak{p}_i \equiv 1$ mod 4 and the corresponding embeddings $\iota_{\infty_j} : k \hookrightarrow \mathbb{R}$ are defined by $\iota_{\infty_1}(a + b\sqrt{p}) = a + b\sqrt{p}$, $\iota_{\infty_2}(a + b\sqrt{p}) = a - b\sqrt{p}$ for $a, b \in \mathbb{Q}$.
In Section \ref{fieldp1p2}, we showed as our first result that the pro-2 Galois group $G_S(2)$ has the following presentation
\begin{equation}\label{intro1}
G_S(2) = \langle x_1, \cdots , x_s \mid x_1^{\mathrm{N}\mathfrak{p}_{1}-1} [ x_1, y_1] = \cdots = x_s^{\mathrm{N}\mathfrak{p}_{s}-1} [ x_s, y_s] =1, x_{\infty_2}^2 = 1 \rangle,
\end{equation}
where the word $x_{i}$ represents a monodromy $\tau_i$ over $\mathfrak{p}_{i}$ in $k_S(2)/k$ and $y_{i}$ denotes the free pro-$2$ word of $x_{1}, \dots, x_{s}$ which represents a Frobenius automorphism $\sigma_i$ over $\mathfrak{p}_{i}$ in $k_S(2)/k$, and $x_{\infty_2}$ denotes the word which represents the monodromy over a real prime $\infty_2$ (see Theorem \ref{maintheorem1}).
We note that the presentation \eqref{intro1} looks similar to that of a link group (\cite{Milnor1957}, \cite{Turaev1979}) by the analogy between a monodromy (resp. Frobenius automorphism) over a prime and a meridian (resp. longitude) of a knot (\cite[Chapter 8]{Morishita2024}).

In Section \ref{deftriplequadraticresiduesymbol}, we define the triple quadratic residue symbol $[\frak{p}_1, \frak{p}_2, \frak{p}_3]$ using group presentation obtained in Section \ref{fieldp1p2}, and we also confirm that this definition is a generalization of the quadratic residue symbol.
To be precise, let us consider the case that $S = \{ \frak{p}_1, \frak{p}_2, \frak{p}_3 \}$ with $\frak{p}_i = (\pi_i)$ such that
{\small \begin{equation*}
\displaystyle{\left( \frac{\pi_i}{\pi_j} \right) = 1} \ {\rm for}\  1 \leq i \neq j \leq 3 \ {\rm and}\  \displaystyle{\left( \frac{\varepsilon}{\pi_i} \right)} = 1 \ {\rm for}\ 1 \leq i \leq 3,
\end{equation*}}
where $\displaystyle{\left(\frac{\cdot}{\cdot} \right)}$ denotes the quadratic residue symbol in $k$.
Then, following the method by Morishita (\cite{Morishita2000}, \cite{Morishita2002}, \cite{Morishita2004}), we obtain the well-defined mod 2 Milnor invariant $\mu_2(123)$ 
by using the presentation \eqref{intro1} and hence the triple quadratic residue symbol 
\begin{equation*}
[\frak{p}_1, \frak{p}_2, \frak{p}_3] := (-1)^{\mu_2(123)}.
\end{equation*}
By the definition of $\mu_2(123)$, the triple symbol $[\frak{p}_1, \frak{p}_2, \frak{p}_3]$ describes the decomposition law of $\frak{p}_3$ in a dihedral extension $K$ of degree 8 over $k$, unramified outside $S$ and the infinite primes with ramification index of $\frak{p}_i$ being 2.
We call such an extension $K$ a {\it R\'edei type $D_8$-extenision} over $k$.
Actually we show that R\'edei type extenision $K$ exists uniquely for $\{\frak{p}_1, \frak{p}_2\}$ (see Theorem \ref{existunique}).

In Section \ref{constructredeiextension}, we construct R\'edei type $D_8$-extenision $K$ under several conditions, following the method by R\'{e}dei (\cite{Amano2014}, \cite{Redei1939}).
To be precise, let $k=\mathbb{Q}(\sqrt{p})$ be a real quadratic field and we assume $p \equiv 5\ {\rm mod}\ 8$.
Let $\frak{p}_1 = (\pi_1), \frak{p}_2 = (\pi_2)$ are prime ideals of $\mathcal{O}_k$ such that
$\pi_i \equiv 1\ {\rm mod}\ 4\mathcal{O}_k$ and $\pi_i > 0\ (1 \leq i \leq 2).$
Then, there is a non-zero $\mathcal{O}_k$ solution $(x, y, z)$ of $x^2 - \pi_1 y^2 - \pi_2 z^2 = 0$ such that ${\rm gcd}(x, y, z) = 1$, $y \in 2\mathcal{O}_k$.
Furthermore, if $x - y \equiv 1\  {\rm mod}\  4\mathcal{O}_k$, then the R\'edei type $D_8$-extenision $K$ is given by
\begin{equation*}
K = k(\sqrt{\pi_1}, \sqrt{\pi_2}, \sqrt{x+y\sqrt{\pi_1}}).
\end{equation*}
We also give some numerical examples of R\'edei type $D_8$-extenisions over real quadratic fields.
Moreover, we give the first example of a {\it Borromean primes} $\frak{p}_1 = (\pi_1), \frak{p}_2 = (\pi_2), \frak{p}_3 = (\pi_3)$ in $\mathbb{Q}(\sqrt{5})$, which satisfies $\left(\frac{\pi_i}{\pi_j} \right) = 1$ ($i \neq j$) and $[\frak{p}_1, \frak{p}_2, \frak{p}_3] = -1$.

Finally, in Section \ref{masseymassey}, following \cite{Morishita2004}, we give a cohomological interpretation of mod 2 Milnor invariant and the triple quadratic residue symbol from the perspective of Massey products.

\subsection*{Notations}
\noindent
For a set $A$, $\#A$ denotes the cardinality of $A$.\\
For a finite algebraic number field $F$,
${\mathcal O}_F$ denotes the ring of integers of $F$,
$H_F$ denotes the ideal class group of $F$,
and $H_F^+$ denotes the narrow ideal class group of $F$.\\
For an ideal $\frak{a}$ of ${\mathcal O}_F$,
${\rm N}\frak{a} := \#({\mathcal O}_k/\frak{a})$.\\
For a group $G$ and $x, y \in G$,
$[x,y] := xyx^{-1}y^{-1}$.\\
For a topological group $G$ and subsets $X, Y$ of $G$,
$[X, Y]$ denotes the closed subgroup of $G$ generated by $[x,y]$ for $x \in X, y \in Y$.

\section{Maximal pro-2 Galois groups of real quadratic fields with restricted ramification}\label{fieldp1p2}

In this section, we study a presentation of the Galois group $G_S(2)$ of the maximal pro-2 extension over a real quadratic field $k$, which is unramified outside a finite set $S$ of primes of $k$.
For this, we apply theorems due to H\"{o}chsmann and Koch (\cite{Koch}) and we assume that the narrow class group of $k$ is trivial and a certain cohomological obstruction $B_S$ vanishes, in order to give a precise minimal generators and their relations.
We also give a condition for $B_S = \{1\}$ in terms of the fundamental unit of $k$.\\

Let $k$ be a real quadratic field and $\mathcal{O}_k$ be the ring of integers of $k$.
We recall that the unit group $\mathcal{O}_k^{\times}$ is given by
$$\mathcal{O}_k^{\times}=\{ \pm \varepsilon^n \mid n \in \mathbb{Z}\},$$
where $\varepsilon$ is the fundamental unit of $k$.
Let $S$ be a finite set of finite primes of $k$,
$$S = \{\frak{p}_1, \cdots, \frak{p}_s\}, \ \#S = s.$$
We assume that $\mathfrak{p}_i$'s are not lying over 2, that is, ${\mathrm N} {\mathfrak p}_i \equiv 1 \ \mathrm{mod} \  2$ for each $i$.
Let $\overline{S}$ be the union
$$\overline{S} = S \cup S_{\mathbb{R}},$$
where $S_{\mathbb{R}} = \{\infty_1, \infty_2\}$ is the set of real primes of $k$.
For $\frak{p}_i \in {S}$, let $k_{\frak{p}_i}$ be the completion of $k$ at $\frak{p}_i$ and let $\pi_i$ be a prime element of $k_{\frak{p}_i}$.
Each primes $\frak{p} \in \overline{S}$ defines an embedding $k$ to $k_{\frak{p}}$, which we denote by $\iota_\frak{p}$.
In particular, for real primes, we denote by $\iota_{\infty_1}$ the embedding into $\mathbb{R}$ determined from the identity map and the conjugate map by $\iota_{\infty_2}$.
Let $k_S(2)$ denote the maximal pro-$2$ extension of $k$, unramified outside $\overline{S}$ and we set $G_S(2) := \mathrm{Gal}(k_S(2)/k)$.

Now we describe the structure of the pro-$2$ group $G_S(2)$ in a certain unobstructed case.
For this, we first recall a result due to Iwasawa on the local Galois group (\cite[Satz 10.2]{Koch}).
We fix an algebraic closure $\overline{k}_{\frak{p}_i}$ of $k_{\frak{p}_i}$ and an embedding $\bar{k} \hookrightarrow \overline{k}_{\frak{p}_i}$.
Let $k_{\frak{p}_i}(2)$ denote the maximal pro-$2$ extension of $k_{\frak{p}_i}$ in $\overline{k}_{\frak{p}_i}$ and $G_{k_{\frak{p}_i}}(2)$ denote the Galois group of $k_{\frak{p}_i}(2)$ over $k_{\frak{p}_i}$.
In this case, the following identity is known:
$$k_{\frak{p}_i}(2) = k_{\frak{p}_i}( \zeta_{2^n}, \sqrt[2^n]{\pi_i} \; \mid \; n \geq 1 ),$$
where $\zeta_{2^n}$ denotes a primitive $2^n$-th root of unity in $\overline{k}$ such that $(\zeta_{2^a})^{2^b} = \zeta_{2^{a-b}}$ for all $a \geq b$.
The local Galois group $G_{k_{\frak{p}_i}}(2)$ is topologically generated by the monodromy $\tau_i$ and (a lift of) the Frobenius automorphism $\sigma_i$ which are defined by
\begin{equation*}
\begin{array}{ll}
\tau_i(\zeta_{2^n}) := \zeta_{2^n}, & \tau_i(\sqrt[2^n]{\pi_i}) := \zeta_{2^n}\sqrt[2^n]{\pi_i},\\
\sigma_i(\zeta_{2^n}) := \zeta_{2^n}^{{\mathrm N}\frak{p}_i},& \sigma_i(\sqrt[2^n]{\pi_i}) := \sqrt[2^n]{\pi_i},
\end{array}
\end{equation*}
which are subject to the relation
\begin{equation}
\label{relationfiniteprime}
\tau_i^{{\mathrm N}\frak{p}_{i} -1}[\tau_i, \sigma_i] = 1.
\end{equation}
For $\frak{p} \in S_{\mathbb{R}}$, we have $k_{\frak{p}} = \mathbb{R}$, $k_{\frak{p}}(2) = \mathbb{C}$, and
\begin{equation}
\label{relationinfinite}
G_{k_\frak{p}}(2) = \langle \ \tau \mid \tau^2=1 \ \rangle.
\end{equation}
For each $\frak{p} \in \overline{S}$, the fixed embedding $\bar{k} \hookrightarrow \overline{k}_{\mathfrak{p}}$ induces an embedding $k_S(2) \hookrightarrow k_{\mathfrak{p}}(2)$, so a prime $\mathfrak{P}_i$ of $k_S(2)$ lying over $\mathfrak{p}_i$ ($1 \leq i \leq s$). 
We denote by the same letters $\tau_i$ and $\sigma_i$ the images of $\tau_i$ and $\sigma_i$ under $\varphi_{\mathfrak{p}}$, respectively, under the homomorphism 
$$\varphi_{\mathfrak{p}} : G_{k_{\mathfrak{p}}}(2) \longrightarrow G_S(2)$$ 
induced by the embedding $k_S(2) \hookrightarrow k_{\mathfrak{p}}(2)$.
Then $\tau_i$ is a topological generator of the inertia group of $\mathfrak{P}_i$ and $\sigma_i$ is a lift of the Frobenius automorphism of the maximal subextension of $k_S(2)/k$ in which $\mathfrak{P}_i$ is unramified. We call simply $\tau_i$  a {\it monodromy over $\mathfrak{p}_i$} in $k_S(2)/k$ and $\sigma_i$ a {\it Frobenius automorphism over $\mathfrak{p}_i$} in $k_S(2)/k$.  
We note that the restriction of $\tau_i$ to the maximal Abelian subextension $k_S(2)^{\rm ab}$ in $k_S(2)/k$, denoted by $\overline{\tau}_i := \tau_i|_{k_S(2)^{\rm ab}}$,
is given by the Artin symbol
$$ \overline{\tau}_i = ( \eta_i, k_S(2)^{\rm ab}/k),$$
for an id\`{e}le $\eta_i$ of $k$ whose $\mathfrak{p}_i$-component is a primitive $(\mathrm{N} \mathfrak{p}_i-1)$-th root of unity in $k_{\mathfrak{p}_i}^{\times}$ and all of the other components are 1.
Now we assume that the $2$-class group $H_k(2)$ of $k$ is trivial.
Then, by class field theory and Burnside's basis theorem (\cite[Satz 4.10]{Koch}), the monodromies $\tau_1, \dots , \tau_s, \tau_{\infty_1}, \tau_{\infty_2}$ generate topologically the global Galois group $G_S(2)$.
However, they may not be a minimal set of generators in general.
In fact, Shafarevich's theorem (\cite[Satz 11.8]{Koch}) tells us that the minimal number $d(G_S(2))$ of generators of $G_S(2)$ is given by
\begin{equation}
\label{dGs}
d(G_S(2)) = s + \dim_{\mathbb{F}_2} B_S,
\end{equation}
where the obstruction $B_S$ is defined as follows: We define the subgroup $V_S$ of $k^{\times}$ by 
$$V_S := \{ a \in k^{\times} \mid a{\mathcal O}_k = \frak{a}^2, \; a \in (k_{\frak{p}}^{\times})^2 \; (\forall \frak{p} \in \overline{S}) \},$$
where $\mathfrak{a}$ is a fractional ideal of ${\mathcal O}_k$, and then $B_S$ is defined by 
$$B_S := V_S/(k^{\times})^2,$$
which is an elementary Abelian $2$-group.
By Shafarevich's formula \eqref{dGs} on $d(G_S(2))$, $2-\dim_{\mathbb{F}_2} B_S$ generators among $\tau_i$'s have to be removed in order to get a minimal generators of $G_S(2)$.

Now we show a necessary and sufficient condition for $B_S$ to vanish by using the fundamental unit $\varepsilon$ of $k$.
For $a \in k^{\times}$ and a finite prime $\frak{p}$ of $k$ such that $2, a$ and $\frak{p}$ are relatively prime, $\left(\frac{a}{\frak{p}}\right)$ denotes the quadratic residue symbol in $k$.
For a real prime $\infty$ of $k$, $\left(\frac{a}{\infty}\right)$ is defined by $1$ (resp. $-1$) if $\iota_\infty(a) > 0$ (resp. if $\iota_\infty(a) < 0$).
\Prp{}{B_S=1}{
{\rm (1)} We assume that the ideal class group $H_k$ of $k$ is trivial.
Then, the following conditions are equivalent.\\
{\rm (i)} $B_S = \{ 1 \}$.\\
{\rm (ii)} There is a prime $\frak{p} \in \overline{S}$ such that $\displaystyle{\left(\frac{\varepsilon}{\frak{p}}\right)} = - 1$.\\
{\rm (2)} We assume that the narrow ideal class group $H_k^+$ of $k$ is trivial.
Then we have $\left(\frac{\varepsilon}{\infty_2}\right) = -1$ and $B_S = \{ 1 \}$.
}

\begin{proof}
{\rm (1)} Suppose that $V_S = (k^\times)^2$.
Since $\varepsilon$ is not square in $k^\times$, $\varepsilon$ is not included in $V_S$.
Since $\varepsilon{\mathcal O}_k = (1)^2$, there exists $\frak{p} \in \overline{S}$ such that $\varepsilon \notin (k_{\frak{p}}^{\times})^2$.
Consequently, we have $\left(\frac{\varepsilon}{\frak{p}}\right) = - 1$ for some prime $\frak{p} \in \overline{S}$.

Conversely, we assume that there is a prime $\frak{p} \in \overline{S}$ such that $\left(\frac{\varepsilon}{\frak{p}}\right) = - 1$.
Choose any $\alpha \in V_S$.
By the definition of $V_S$, there is an ideal $\frak{a}$ satisfying $\alpha{\mathcal O}_k = \frak{a}^2$.
Since the class number of $k$ is one, $\frak{a}$ is generated by one element $a \in k^{\times}$.
So, $\alpha$ has the form $\alpha = u \cdot a^2$ ($u \in \mathcal{O}_{k}^{\times}$).
Because of the structure of ${\mathcal O}_k^{\times}/({\mathcal O}_k^{\times})^2$, we have $u = u' \cdot v^2$ ($u' \in \{1, -1, \varepsilon, -\varepsilon\}$, $v \in {\mathcal O}_k^{\times}$).
Then replacing $v \cdot a$ by $a$, $\alpha$ has the form $u' \cdot a^2$ ($u' \in \{1, -1, \varepsilon, -\varepsilon\}$).
So we obtain $u' = \alpha \cdot (a^{-1})^2 \in V_S$.
On the other hand, $-1,-\varepsilon \notin V_S$ because $\overline{S}$ has a real prime $\infty_1$.
In addition, $\varepsilon$ is not included in $V_S$ by the assumption.
Consequently, we have $u' = 1$ and $\alpha = a^2 \in (k^{\times})^2$.\\
{\rm (2)} Since the narrow ideal class group $H_k^+$ of $k$ is trivial, ${\rm N}_{k/\mathbb{Q}}(\varepsilon) = \varepsilon \cdot \iota_{\infty_2} (\varepsilon) = -1$.
Therefore we have $\iota_{\infty_2} (\varepsilon) = -1$ and $\left(\frac{\varepsilon}{\infty_2}\right) = -1$.
By {\rm (1)}, we have $B_S = \{ 1 \}$.
\end{proof}

By using Proposition \ref{B_S=1}, the minimal set of generators of $G_S(2)$ can be determined as follows.
Let $\Phi(G_S(2))$ be the Frattini subgroup of $G_S(2)$, namely,
$\Phi(G_S(2)) = G_S(2)^2[G_S(2), G_S(2)].$

\Prp{}{structure of G_S}{Notations being as above, we assume that the ideal class group $H_k$ of $k$ is trivial, and there is a prime $\frak{p} \in \overline{S} \setminus \{\infty_1\}$ such that $\left(\frac{\varepsilon}{\frak{p}}\right) = - 1$.
Then $G_{S}(2)$ is topologically generated by the set of minimal generators $\{\tau_1, \dots, \tau_{s}, \tau_{\infty_2}\}$ $\setminus \{\tau\}$, where $\tau$ is a monodromy corresponding to $\frak{p}$.
If the narrow ideal class group $H_k^+$ of k is trivial, then we can take $\tau = \tau_{\infty_2}$.}
\begin{proof}
By Burnside's basis theorem, it is sufficient to show that $\{\tau_1, \dots, \tau_{s}, \tau_{\infty_2}\} \setminus \{\tau\}$ generates $G_S(2)/\Phi(G_S(2))$.
For a prime $\frak{p}$ of $k$, let ${\mathcal O}_{\frak{p}}$ be the ring of $\frak{p}$-adic integers of $k_{\frak{p}}$.
(For the infinite prime $\frak{p}$, we set ${\mathcal O}_{\frak{p}}:= k_{\frak{p}}$.)
Let $J_k$ be the id\`{e}le group of $k$.
Let $U_k$ be the subgroup of $J_k$ consisting of unit id\`{e}les of $k$, $U_k := \prod_{\frak{p}} {\mathcal O}_{\frak{p}}^{\times}$, and let $U_{S}$ denote the subgroup of $U_k$ whose $\frak{p}$-component is $1$ for any $\frak{p} \in \overline{S}$.
Employing class field theory, we have the canonical isomorphism of $\mathbb{F}_2$-vector spaces 
$$J_k/U_{S}J_k^2 k^{\times} \; \simeq \; G_S(2)/\Phi(G_S(2))$$
 by the Artin symbol.

Since the ideal class group $H_k = J_k/U_kk^{\times}$ is trivial and $B_S = \{1\}$ by Proposition \ref{B_S=1}, we have the following exact sequence of $\mathbb{F}_2$-vector spaces (cf. \cite[(11.11)]{Koch}).
\begin{equation*}
1 \longrightarrow {\mathcal O}_k^{\times}/({\mathcal O}_k^{\times})^2 \stackrel{\delta}{\longrightarrow} \prod_{\frak{p} \in \overline{S}} {\mathcal O}_{\frak{p}_i}^{\times}/({\mathcal O}_{\frak{p}_i}^{\times})^2 \stackrel{\iota}{\longrightarrow} J_k/U_{S}J_k^2 k^{\times} \longrightarrow 1,
\end{equation*}
where $\delta$ is the diagonal map and $\iota$ is induced by the natural inclusion $\prod_{\frak{p} \in \overline{S}} {\mathcal O}_{\frak{p}_i}^{\times} \hookrightarrow J_k$.
By local and global class field theory, ${\mathcal O}_{\frak{p}_i}^{\times}/({\mathcal O}_{\frak{p}_i}^{\times})^2 \simeq \langle \overline{\tau_i} \rangle$ and $J_k/U_{S}J_k^2 k^{\times} \simeq G_S(2)/\Phi(G_S(2))$.
So we obtain the following exact sequence.
\begin{equation*}
1 \rightarrow {\mathcal O}_k^{\times}/({\mathcal O}_k^{\times})^2  \stackrel{\delta}{\longrightarrow} \langle \overline{\tau_1} \rangle \times \cdots \times \langle \overline{\tau_{s}} \rangle \times \langle \overline{\tau_{\infty_1}} \rangle \times \langle \overline{\tau_{\infty_2}} \rangle \stackrel{\iota}{\longrightarrow} G_S(2)/\Phi(G_S(2)) \rightarrow 1.
\end{equation*}
Therefore, $\{\iota (\overline{\tau_1}), \cdots \iota (\overline{\tau_{s}}), \iota (\overline{\tau_{\infty_1}}), \iota (\overline{\tau_{\infty_2}})\}$ forms a system of generators of \\$G_S(2)/\Phi(G_S(2))$.
On the other hand, Since ${\mathcal O}_k^{\times}/({\mathcal O}_k^{\times})^2 = \{1, -1, \varepsilon, -\varepsilon\} = \langle -1, \varepsilon \rangle$, $\{\delta(-1), \delta(\varepsilon)\}$ is a basis of ${\rm Im}(\delta)$.
Since ${\rm Ker}(\iota) = {\rm Im}(\delta)$, $\{\delta(-1), \delta(\varepsilon)\}$ is a basis of ${\rm Ker}(\iota)$.
By expressing $\delta(-1)$ and $\delta(\varepsilon)$ as $\delta(-1) = \overline{\tau_1}^{a_1} \cdots \overline{\tau_{\infty_2}}^{a_{\infty_2}}$, $\delta(\varepsilon) = \overline{\tau_1}^{b_1} \cdots \overline{\tau_{\infty_2}}^{b_{\infty_2}}$ using the basis of $\langle \overline{\tau_1} \rangle \times \cdots \times \langle \overline{\tau_{\infty_2}} \rangle$, we obtain the relations
\begin{equation*}
\begin{cases}
\overline{\tau_1}^{a_1} \cdots \overline{\tau_{\infty_2}}^{a_{\infty_2}} = 1 & \text{(\ref{structure of G_S}.1)}\\
\overline{\tau_1}^{b_1} \cdots \overline{\tau_{\infty_2}}^{b_{\infty_2}}=1 & \text{(\ref{structure of G_S}.2)} 
\end{cases}
\end{equation*}
in $G_S(2)/\Phi(G_S(2))$.
Here, for $i, j = 1, \cdots, s, \infty_1, \infty_2$, we define
\begin{equation*}
a_i = \begin{cases}
1 & \delta(-1) \notin ({\mathcal O}_{\frak{p}_i}^{\times})^2 \\
0 & \text{otherwise,}
\end{cases}
\ b_j = \begin{cases}
1 & \delta(\varepsilon) \notin ({\mathcal O}_{\frak{p}_j}^{\times})^2 \\
0 & \text{otherwise.}
\end{cases}
\end{equation*}
From relation (\ref{structure of G_S}.1), since $\delta(-1) = -1 \notin ({\mathcal O}_{\infty_1}^{\times})^2 = (\mathbb{R}^{\times})^2$ we see that it is possible to omit $\overline{\tau_{\infty_1}}$ from the system of generators of $G_S(2)/\Phi(G_S(2))$.
Moreover, by defenition of $b_j$ and relation (\ref{structure of G_S}.2) another omissible element is $\overline{\tau}$ corresponding to $\frak{p} \in \overline{S}$ such that $\left(\frac{\varepsilon}{\frak{p}}\right) = - 1$.

In addition, if the narrow ideal class group of $k$ is trivial, the same argument as above shows that we can use the relation (\ref{structure of G_S}.1) (resp. (\ref{structure of G_S}.2)) to omit $\tau_{\infty_1}$ (resp. $\tau_{\infty_2}$) from the system of generators.
Note that the image of $\varepsilon$ into ${\mathcal O}_{\infty_2}^{\times}$ by $\delta$ is the conjugate of $\varepsilon$, and the conjugate of $\varepsilon$ is a negative number when the narrow ideal class group is trivial.
\end{proof}

By Proposition \ref{B_S=1} and Proposition \ref{structure of G_S}, the following holds.

\Thm{}{maintheorem1}{Notations being as above, we assume that the narrow ideal class group $H_k^+$ of $k$ is trivial.
Then the pro-$2$ group $G_S(2)$ has the following presentation:
\[
G_S(2) = \langle x_1, \cdots , x_s \mid x_1^{\mathrm{N}\mathfrak{p}_{1}-1} [ x_1, y_1] = \cdots = x_s^{\mathrm{N}\mathfrak{p}_{s}-1} [ x_s, y_s] =1, x_{\infty_2}^2 = 1 \rangle,
\]
where $x_{i}$ denotes the word which represents a monodromy $\tau_i$ over $\mathfrak{p}_{i}$ in $k_S(2)/k$ and $y_{i}$ denotes the free pro-$2$ word of $x_{1}, \dots, x_{s}$ which represents a Frobenius automorphism $\sigma_i$ over $\mathfrak{p}_{i}$ in $k_S(2)/k$, and $x_{\infty_2}$ denotes the word which represents the monodromy over a real prime $\infty_2$.
}
\begin{proof}
Since the narrow ideal class group $H_k^+$ of $k$ is trivial, $B_S = \{ 1 \}$ by Proposition \ref{B_S=1} and $G_{S}(2)$ is topologically generated by the set of minimal generators $\{\tau_1, \dots, \tau_{s}\}$ by Proposition \ref{structure of G_S}.
Since $B_S = \{ 1 \}$, according to \cite[Satz 11.3]{Koch}, the localization map
\begin{equation*}
H^2(G_S(2), \mathbb{F}_2) \longrightarrow \prod_{\frak{p} \in \overline{S}} H^2(G_{k_{\frak{p}}}(2), \mathbb{F}_2)
\end{equation*}
is injective.
Moreover, \cite[Satz 11.4]{Koch} tells us that the map
\begin{equation*}
H^2(G_S(2), \mathbb{F}_2) \longrightarrow \prod_{\frak{p} \in \overline{S} \setminus \{\infty_1\}} H^2(G_{k_{\frak{p}}}(2), \mathbb{F}_2)
\end{equation*}
is also injective.
Therefore, by \cite[Satz 6.14]{Koch}, the relation system of $G_S(2)$ is generated by the local relation system of $G_{k_{\frak{p}}}(2)$ ($\frak{p} \in \overline{S} \setminus \{\infty_1\}$).
By \eqref{relationfiniteprime} and \eqref{relationinfinite}, we have presentation of $G_S(2)$.
\end{proof}

The following Proposition \ref{xinfty in frattini} will be used in Section \ref{deftriplequadraticresiduesymbol}.

\Prp{}{xinfty in frattini}{Notations being as above, we assume that the narrow ideal class group $H_k^+$ of $k$ is trivial.
Then, the following conditions are equivalent.\\
{\rm (i)} $x_{\infty_2} \in \Phi(G_S(2))$.\\
{\rm (ii)} $\displaystyle{\left(\frac{\varepsilon}{\frak{p}_i}\right)} = 1$ for all finite primes $\frak{p}_i \in S$.
}
\begin{proof}
By the proof of Proposition \ref{structure of G_S}, we have relations
$$\begin{cases}
\overline{\tau_1}^{a_1} \cdots \overline{\tau_{\infty_2}}^{a_{\infty_2}} = 1, \\
\overline{\tau_1}^{b_1} \cdots \overline{\tau_{\infty_2}}^{b_{\infty_2}}=1
\end{cases}$$
in $G_S(2)/\Phi(G_S(2))$.
Here, for $i, j = 1, \cdots, s, \infty_1, \infty_2$, we define
$$a_i = \begin{cases}
1 & \delta(-1) \notin ({\mathcal O}_{\frak{p}_i}^{\times})^2 \\
0 & \text{otherwise,}
\end{cases}
\ b_j = \begin{cases}
1 & \delta(\varepsilon) \notin ({\mathcal O}_{\frak{p}_j}^{\times})^2 \\
0 & \text{otherwise.}
\end{cases}$$
Therefore, $x_{\infty_2} \in \Phi(G_S(2))$ is equivalent to that one of the following conditions hold:\\
(iii) for $i = 1, \cdots, s, \infty_1$, \[-1 \in ({\mathcal O}_{\frak{p}_i}^{\times})^2 \ {\rm and}\  -1 \notin ({\mathcal O}_{\frak{p}_{\infty_2}}^{\times})^2,\]\\
(iv) for $j = 1, \cdots, s, \infty_1$, \[\varepsilon \in ({\mathcal O}_{\frak{p}_j}^{\times})^2 \ {\rm and}\  \varepsilon \notin ({\mathcal O}_{\frak{p}_{\infty_2}}^{\times})^2.\]
Since $-1 \notin ({\mathcal O}_{\frak{p}_{\infty_1}}^{\times})^2$, (iii) cannot hold.
An equivalent condition for (iv) is that \[\left(\frac{\varepsilon}{\frak{p}_i}\right) = 1\]
for all finite primes $\frak{p}_i \in S$.
\end{proof}

\section{Mod 2 arithmetic Milnor invariants and triple quadratic residue symbols in real quadratic fields}\label{deftriplequadraticresiduesymbol}
In this section, following Morishita's work (\cite{Morishita2000}, \cite{Morishita2002}, \cite{Morishita2004}), we introduce mod 2 arithmetic Milnor invariants and triple quadratic residue symbols for certain primes of a real quadratic field, by using the presentation in Theorem \ref{maintheorem1} in Section \ref{fieldp1p2}.\\

Now, let us be back to the arithmetic situation in Section \ref{fieldp1p2} and keep the same notations and assumptions there.
Let $k$ be a real quadratic field.
In the following, we assume that the narrow ideal class group $H_k^+$ of $k$ is trivial.
Let $S = \{\mathfrak{p}_1, \cdots, \mathfrak{p}_s\}$ be a set of $s$ distinct finite primes of $k$, where $\mathfrak{p}_i$'s are not lying over 2, and let $\overline{S} := S \cup S_{\mathbb{R}}$, when $S_{\mathbb{R}}$ denotes the set of real primes of $k$.
We fix a prime element $\pi_i$ of $\frak{p}_i$, namely, $\frak{p}_i = (\pi_i)$.
By Theorem \ref{maintheorem1}, the pro-2 Galois group $G_S(2)$ has the following presentation:
\[
G_S(2) = \langle x_1, \cdots , x_s \mid x_1^{\mathrm{N}\mathfrak{p}_{1}-1} [ x_1, y_1] = \cdots = x_s^{\mathrm{N}\mathfrak{p}_s-1} [ x_s, y_s] =1, x_{\infty_2}^2 = 1 \rangle,
\]
where $x_{i}$ denotes the word which represents a monodromy over $\mathfrak{p}_{i}$ in $k_S(2)/k$ and $y_{i}$ denotes the free pro-$2$ word of $x_{1}, \dots, x_{s}$ which represents a Frobenius automorphism over $\mathfrak{p}_{i}$ in $k_S(2)/k$, and $x_{\infty_2}$ denotes the word which represents a monodromy over a real prime $\infty_2$.

Let $F_S$ denote the free pro-2 group on the words $x_1, \dots, x_s$.
Let $\mathbb{F}_2[[F_S]]$ be the complete group algebra of $F_S$ over $\mathbb{F}_2$ and let $\varepsilon_{\mathbb{F}_2[[F_S]]} : \mathbb{F}_2[[F_S]] \longrightarrow \mathbb{F}_2$ be the augmentation homomorphism with the augmentation ideal $I_{\mathbb{F}_2[[F_S]]} := {\rm Ker}(\varepsilon_{\mathbb{F}_2[[F_S]]})$.
Let $\mathbb{F}_2 \langle \langle X_1, \dots, X_s \rangle \rangle$ denote the formal power series algebra over $\mathbb{F}_2$ in non-commutative variables $X_1, \dots , X_s$.
For a multiple index $I=(i_1 \cdots i_n), 1\leq i_1,\dots ,i_n \leq s$, we set $X_I = X_{i_1} \cdots X_{i_n}$.
Let 
$$ M : \mathbb{F}_2[[F_S]] \; \stackrel{\sim}{\longrightarrow} \; \mathbb{F}_2\langle\langle X_1,\dots ,X_s \rangle\rangle,$$
be {\it the {\rm mod} $2$ Magnus isomorphism} defined by
$$ M(x_i) = 1+X_i,\enspace M(x_i^{-1}) = 1-X_i+X_i^2- \cdots .$$
Then, for $f \in F_S$, $M(f)$ has the form (called {\it {\rm Mod} $2$ Magnus expansion})
$$M(f) = 1 + \sum_{|I| \geq 1} \mu_2(I;f) X_I,$$
where $\mu_2(I;f) \in \mathbb{F}_2$ are the mod $2$ Magnus coefficients.
In terms of the pro-2 Fox free derivatives $\frac{\partial}{\partial x_i} : \mathbb{Z}_2[[F_S]] \longrightarrow \mathbb{Z}_2[[F_S]]$ over $\mathbb{Z}_2$, we can obtain for $I=(i_1 \cdots i_n)$ $$\mu_2(I;f) = \varepsilon_{{\mathbb{F}_2[[F_S]]}}\left(\frac{\partial^n f}{\partial x_{i_1} \cdots \partial x_{i_n}} \ {\rm mod}\  2 \right),$$
where mod $2$ : $\mathbb{Z}_2[[F_S]] \longrightarrow \mathbb{F}_2[[F_S]]$ is the reduction mod $2$.

Here, we recall basic properties of mod $2$ Magnus coefficients.

\Lem{}{magnuscoeffproperty}{{\rm (i)} For $\alpha, \beta \in \mathbb{F}_2[[F_S]]$ and a multi-index $I$,
$$\mu_2(I;\alpha\beta) = \sum_{I=JK} \mu_2(J;\alpha)\mu_2(K;\beta),$$
where the sum ranges over all pairs $(J,K)$ of multi-indices such that $JK = I$.\\
{\rm (ii)} For $f \in F_S$ and $q \geq 2$, we have
$$\mu_2(I;f) = 0 \ {\rm for}\  |I| < q \iff f \in {F_S}^{(q)},$$
where $F_S^{(q)} := \{ x \in F_S \mid  g - 1 \in I_{\mathbb{F}_2[[F_S]]}^q\}$, the $q$-th term of the mod $2$ Zassenhaus filtration of $F_S$.}

We note that the set $\{F_S^{(q)}\}$ forms a fundamental system of neighbourhoods of 1 in $F_S$ and that $F_S^{(q)}/F_S^{(q+1)}$ is central in $F_S/F_S^{(q+1)}$ and $(F_S^{(q)})^2 \subset F_S^{(2q)}$.

Now let $I = (i_1 \cdots i_n)$ be a multi-index with $1\leq i_1, \dots , i_n \leq s$.
Let $I' := (i_1 \cdots i_{n-1})$.
Then {\it the mod $2$ arithmetic Milnor number} $\mu_2(I)$ for $I$ is defined by
$$ \mu_2(I) := \mu_2(I'; y_{i_n}).$$
Here we set $\mu_2(I) := 0$ if $|I| = 1$.
The mod $2$ arithmetic Milnor numbers $\mu_2(I)$ are not invariants determined by the Galois group $G_S(2)$.
However, we can show in the following that $\mu_2(I)$ are invariant of $G_S(2)$ if $|I| = 2$ or $3$.
\Prp{}{millnor residue}{Notations being as above, let $i, j$ be indices between $1$ and $s$.
When $i=j$, we have $\mu_2(i i) = 0$.
When $i \neq j$, we have
$$(-1)^{\mu_2(i j)} = \left(\frac{\pi_i}{\pi_j}\right).$$
}
\begin{proof}
By definition of the quadratic residue symbol,
\begin{align*}
\left(\frac{\pi_i}{\pi_j}\right) &= \frac{(\pi_j, k_{\frak{p}_j}(\sqrt{\pi_i})/k_{\frak{p}_j})(\sqrt{\pi_i})}{\sqrt{\pi_i}}\\
&= \frac{y_j(\sqrt{\pi_i})}{\sqrt{\pi_i}}.
\end{align*}
where $(\quad ,k_{\frak{p}_j}(\sqrt{\pi_i})/k_{\frak{p}_j}) : k_{\frak{p}_j}^{\times} \longrightarrow {\rm Gal}(k_{\frak{p}_j}(\sqrt{\pi_i})/k_{\frak{p}_j})$ is the norm residue symbol.
On the other hand, we consider the mod 2 Magnus isomorphism of $y_j$,
\begin{equation*}
M(y_j) = 1 + \sum_i \mu_2(i j) X_i + {\rm (terms \ of \ higher \ degree)}.
\end{equation*}
Therefore $y_j \equiv \prod_i x_i^{\mu_2(i j)} \ {\rm mod} \ [G_S(2),G_S(2)]$.
Hence
\begin{align*}
\frac{y_j(\sqrt{\pi_i})}{\sqrt{\pi_i}} &= \frac{x_i^{\mu_2(i j)} (\sqrt{\pi_i})}{\sqrt{\pi_i}}\\
&= (-1)^{\mu_2(i j)}.
\end{align*}
\end{proof}
Next, we show that mod $2$ arithmetic triple Milnor numbers are invariants determined by $G_S(2)$.
\Thm{}{invarianttheorem}{Notations being as above, we assume that $\left(\frac{\varepsilon}{\frak{p}_i}\right) = 1$ for all finite primes $\frak{p}_i \in S$, $\mu_2(J) = 0$ for any multi-index $J$ of length $\leq 2$ ,and ${\rm N}\frak{p}_i \equiv 1$ {\rm mod} $4$ for all $\frak{p}_i \in S$.
Then, for a multi-index $I$ satisfying $|I| = 3$, $\mu_2(I)$ is an invariant depending on $G_S(2)$.}
\begin{proof}
Let $I := (i_1, i_2, i_3)$.
We must show that $\mu_2(I)$ is independent of the choices of a monodromy over $\frak{p}_i$ and an extension of the Frobenius automorphism over $\frak{p}_i$.
It suffices to show the following:\\
($\mathrm{i}$) $\mu_2(I)$ is not changed if $y_{i}$ is replaced by a conjecture.\\
($\mathrm{ii}$) $\mu_2(I)$ is not changed if $x_{i}$ is replaced by a conjecture.\\
($\mathrm{iii}$) $\mu_2(I)$ is not changed if $y_{i_3}$ is multiplied by a conjugate of $x_i^{{\rm N}\frak{p}_{i}-1} [x_i, y_i]$ and $x_{\infty_2}$.\\
Let $I' := (i_1 , i_2)$.\\
{\bf Proof of ($\mathrm{i}$)}: Write $M(y_i) = 1 + \omega_i$ for the Mod 2 Magnus expansion of $y_i$, $i = 1, 2, 3$.
By the assumption, Lemma \ref{magnuscoeffproperty} (ii) and Proposition \ref{millnor residue}, the degree of $\omega_i \geq 2$.
Then, we have
\begin{align*}
 M(x_jy_ix_j^{-1}) &= (1 + X_j)(1 + \omega_i)(1-X_j+X_j^2- \cdots ) \\
                            &= 1+ \omega_i + ({\rm terms \ involving} \ X_j\omega_i \ {\rm or}\  \omega_iX_j) \\
                            &\equiv M(y_i) + ({\rm terms \ of \ degree} \geq 3).
\end{align*}
Since $F_S$ is topologically generated by $x_j$'s, this proves (i).\\
{\bf Proof of ($\mathrm{ii}$)}: Suppose that $x_i$ is replaced by $x_i^{*} = x_jx_ix_j^{-1}$.
As $x_i = x_j^{-1}x_i^{*}x_j$, we have $1 + X_i = (1-X_j+X_j^2- \cdots ) (1+X_i^{*}) (1+X_j)$.
Therefore, $X_i = X_i^{*} + ({\rm terms \  containing}\  X_jX_i^{*}\  {\rm or}\  X_i^{*}X_j)$.
Each time $X_i$ occurs in the Magnus expansion $M(y_{i})$, it is to be replaced by this last expansion.
The terms in the bracket give rise to terms of degree $\geq 3$ in the new Magnus expansion of $y_i$ using $X_i^*$.
Hence, the new expansion has the same coefficient in degree $2$ as the old.
This completes the proof of (ii).\\
{\bf Proof of ($\mathrm{iii}$)}: The identity $M(x_jy_j-y_jx_j) = X_j \omega_j - \omega_j X_j$ implies
\begin{align*}
M([x_j, y_j]) &= M(1 + (x_jy_j - y_jx_j) x_j^{-1}y_j^{-1})\\
              &= 1+ (X_j \omega_j - \omega_j X_j)(1-X_j+X_j^2- \cdots )M(y_j^{-1})\\
              &= 1 \ {\rm mod} \ M(I_{\mathbb{F}_2[[F_S]]}^3).
\end{align*}
In addition, we have $M(x_j^{{\rm N}\frak{p}_j-1}) = (1+ X_j)^{{\rm N}\frak{p}_j-1} \equiv 1 \ {\rm mod} \ M(I_{\mathbb{F}_2[[F_S]]}^3)$.
Hence, we have
\begin{equation*}
x_j^{{\rm N}\frak{p}_j-1} [x_j, y_j] \equiv 1 \ {\rm mod} \ I_{\mathbb{F}_2[[F_S]]}^3.
\end{equation*}
Moreover, by Proposition \ref{xinfty in frattini}, $x_{\infty_2} \in \Phi(G_S(2)) = G_S(2)^{(2)} \subset F_S^{(2)}$ and $(F_S^{(2)})^2 \subset F_S^{(4)}$, $x_{\infty_2}^2 \in F_S^{(4)}$.
Hence we have
\begin{equation*}
M(x_{\infty_2}^2) = 1 + ({\rm terms \ of \ degree} \geq 4),
\end{equation*}
which leads to the assertion.
\end{proof}

By Theorem \ref{invarianttheorem}, we call $\mu_2(I)$ the mod 2 {\it arithmetic triple Milnor invariant} of $S$ for the multi-index $|I| = 3$, when $\left(\frac{\varepsilon}{\frak{p}_i}\right) = 1$ for $\frak{p}_i \in S$, $\mu_2(J) = 0$ for any multi-index $J$ of length $\leq 2$ and ${\rm N}\frak{p}_i \equiv 1$ {\rm mod} $4$ for $\frak{p}_i \in S$.

A meaning of the mod 2 arithmetic triple Milnor invariants is given as follows.
Let $\overline{S} = S \cup S_{\mathbb{R}} = \{ \frak{p}_1, \frak{p}_2 , \frak{p}_3\} \cup \{\infty_1,\infty_2\}$ be a finite set of primes of $k$, where ${\rm N}\frak{p}_i \equiv 1 \ {\rm mod} \ 4$ ($1 \leq i \leq 3$) and $\infty_1, \infty_2$ are real primes.
We assume that $\mu_2(ij) = 0$ for any multi-index of length $2$.
By Proposition \ref{millnor residue}, it means $\left(\frac{\pi_i}{\pi_j}\right)=1$ for $i \neq j$.
Let $N_3(\mathbb{F}_2)$ denote the group of $3$ by $3$ upper-triangular unipotent matrices over $\mathbb{F}_2$.
We note that $N_3(\mathbb{F}_2)$ is isomorphic to the dihedral group $D_8 = \langle s, t \; | \; s^2 = t^4 = 1, sts^{-1} = t^{-1} \rangle$ of order $8$ by the correspondence
\begin{equation*}
s \mapsto \left( \begin{array}{ccc} 1 & 0 & 0 \\ 0 & 1 & 1 \\ 0 & 0 & 1 \end{array} \right), \;\; t \mapsto \left( \begin{array}{ccc} 1 & 1 & 0 \\ 0 & 1 & 1 \\ 0 & 0 & 1 \end{array} \right).
\end{equation*}
We define the map
\begin{equation*}
\rho : F_{S} \longrightarrow N_3(\mathbb{F}_2)
\end{equation*}
by
\begin{equation*}
\rho(f) := \left(
\begin{array}{ccc}
	 1 & \mu_2((1);f) & \mu_2((1 2); f)\\
	 0 & 1 & \mu_2((2); f)\\
	 0 & 0 & 1 \\
\end{array}
\right)
\end{equation*}
for $f \in F_S$.
By Lemma \ref{magnuscoeffproperty} (i), we see $\rho$ is a group homomorphism.
Let $S' = \{\frak{p}_1, \frak{p}_2\}$ and let $\overline{S'} = S' \cup S'_{\mathbb{R}}$ be the subset of $S$ and $\overline{S}$.

\Thm{}{presentation}{Notations being as above, we assume that $\left(\frac{\varepsilon}{\frak{p}_i}\right) = 1$ for $\frak{p}_i \in S$.\\
{\rm (1)} The homomorphism $\rho$ is surjective and factors through the Galois group $G_S(2)$ and $G_{S'}(2)$.\\
{\rm (2)} Let $K_\rho$ be the extension over $k$ corresponding to {\rm Ker}{\rm($\rho$)}.
Then $K_\rho$ is a Galois extension of $k$ unramified outside $\overline{S'}$ with Galois group ${\rm Gal}(K_\rho/k) = N_3(\mathbb{F}_2)$, and each ramification index over $\frak{p}_i$ is $2$.\\
{\rm (3)} Assume that $\left(\frac{\pi_i}{\pi_j}\right)=1$ for $i \neq j$.
For a Frobenius automorphism $\sigma_{3}$ over $\frak{p}_{3}$, we have
\begin{equation*}
\rho(\sigma_{3}) = \left(
\begin{array}{ccc}
	 1 &  0 & \mu_2(123)\\
	 0 & 1 &  0\\
	 0 & 0 & 1 \\
\end{array}
\right).
\end{equation*}
In particular, we have
\begin{equation*}
\mu_2(1 2 3) = 0 \Longleftrightarrow \frak{p}_{3}\  \mbox{is completely decomposed in} \  K_\rho/k.
\end{equation*}
}           
\begin{proof}
(1) Since
\begin{equation*}
\begin{array}{cc}
 \rho(x_1) = \left(
\begin{array}{ccc}
	 1 & 1 & 0 \\
	 0 & 1 & 0 \\
	 0 & 0 & 1 \\
\end{array}
\right), & \rho(x_2) = \left(
\begin{array}{ccc}
	 1 & 0 & 0 \\
	 0 & 1 & 1 \\
	 0 & 0 & 1 \\
\end{array}
\right)
\end{array}
\end{equation*}
and so $\rho(x_1), \rho(x_2)$ generate $N_3(\mathbb{F}_2)$, $\rho$ is surjective.\\
Next, by presentation of $G_S(2)$ and $G_{S'}(2)$, we need to show
\begin{equation*}
\rho(x_i^{\mathrm{N}\frak{p}_i-1}[x_i,y_i]) = E_3 (1 \leq i \leq 3) \ {\rm and}\  \rho(x_{\infty_2}^2) = E_3,
\end{equation*}
where $E_3$ is the identity matrix of size 3.
Since we can prove \[\mu_2((1); x_i^{\mathrm{N}\frak{p}_i-1}[x_i,y_i]) = \mu_2((2); x_i^{\mathrm{N}\frak{p}_i-1}[x_i,y_i]) = \mu_2((12); x_i^{\mathrm{N}\frak{p}_i-1}[x_i,y_i])= 0\] and \[\mu_2((1); x_{\infty_2}^2) = \mu_2((2); x_{\infty_2}^2) = \mu_2((12); x_{\infty_2}^2)= 0\] as in the proof of Theorem \ref{invarianttheorem}, we are done.\\
(2) By (1), we obtain Gal($K_\rho/k$) = $N_3(\mathbb{F}_2)$.
Since $\rho(x_3) = E_3$, $K_\rho/k$ is unramified outside $\overline{S'}$.
The ramification index over $\frak{p}_i$ is the order of $\rho(x_i)$, which is 2.\\
(3) For $J$ with $|J| = 1$, $\mu_2(J;y_3) = 0$ by definition and the assumption.
Hence the assertion follows.
\end{proof}
\Def{}{}{
Let $\frak{p}_1, \frak{p}_2$ be distinct finite primes of $k$ such that ${\rm N}\frak{p}_i \equiv 1 \ {\rm mod}\  4$.
We call an extension $K$ of $k$ a {\it R\'edei type $D_8$-extension of $k$ for $\{\frak{p}_1, \frak{p}_2\}$} if $K/k$ is a Galois extension of $k$ which is unramified outside $\{\frak{p}_1, \frak{p}_2, \infty_1, \infty_2\}$ with ramification index of each $\frak{p}_i$ being 2 and whose Galois group is isomorphic to the group $D_8 \cong N_3(\mathbb{F}_2)$.
}
\Thm{}{existunique}{Notations being as above,
let $\frak{p}_1 = (\pi_1), \frak{p}_2 = (\pi_2)$ be distinct finite primes of $k$ such that ${\rm N}\frak{p}_i = 1$ {\rm mod} $4$.
Suppose that $\left(\frac{\pi_1}{\pi_2}\right) = \left(\frac{\pi_2}{\pi_1}\right) = 1$.
$\iota_{\infty_2}(\pi_1) > 0$, and $\iota_{\infty_2}(\pi_2) > 0$\\
{\rm (1) (Uniqueness)} The R\'edei type $D_8$-extension of $k$ for $\{\frak{p}_1, \frak{p}_2\}$ is unique if it exists.\\
{\rm (2) (Existence)} We assume further $\left(\frac{\varepsilon}{\frak{p}_i}\right) = 1$ for $\frak{p}_i$ for $i=1,2$.
Then the R\'edei type $D_8$-extension of $k$ for $\{\frak{p}_1, \frak{p}_2\}$ exists.}
\begin{proof}
(1) This is proved in a manner similar to the proof of \cite[Theorem 4.1]{AMM}, \cite[Proposition 3.7]{Mizusawa}.
Let $K/k$ be a R\'edei type $D_8$-extension of $k$ for $\{\frak{p}_1, \frak{p}_2\}$ and $S' = \{\frak{p}_1 = (\pi_1), \frak{p}_2 = (\pi_2)\}$.
We fix a presentation of $G_{S'}(2)$:
\begin{equation*}
G_{S'}(2) = F/R,
\end{equation*}
where $F$ is the free pro-2 group on the words $x_1, x_2$ and $R$ is a normal subgroup of $F$ generated by $x_1^{\mathrm{N}\frak{p}_1-1}[x_1,y_1]$, $x_2^{\mathrm{N}\frak{p}_2-1}[x_2,y_2]$, $x_{\infty_2}^2$.
Let $N$ be the normal subgroup of $F$ generated by $x_1^2, x_2^2$.
Since any element of $F^2$ can be expressed as a product of $x_1^2$ and $x_2^2$ modulo $[F, F]$, $N[F, F] = F^2[F, F]$.
By definition of $R$, $R$ is subgroup of $N[F, F] = F^2[F, F]$.
By Proof of Proposition \ref{millnor residue}, our assumption $\left(\frac{\pi_1}{\pi_2}\right) = \left(\frac{\pi_2}{\pi_1}\right) = 1$ implies that $y_1, y_2 \in [F, F] \subset N[F,F]$, so $[x_1, y_1], [x_2, y_2] \in N[[F, F], F]$.
In addition, since $N[F, F]$ corresponds to $k(\sqrt{\pi_1}, \sqrt{\pi_2})$, we have $x_{\infty_2} \in N[F, F]$ by our assumption $\iota_{\infty_2}(\pi_1) > 0$ and $\iota_{\infty_2}(\pi_2) > 0$.
In otherwords, $x_{\infty_2}^2 \in (N[F, F])^2 \subset N[[F, F], F]$.
Since
\begin{equation*}
[x_1, x_2]^2 \equiv [x_1^2,x_2] \equiv 1\ {\rm mod}\ N[[F, F], F],
\end{equation*}
$N[F, F]/N[[F, F], F]$ is a cyclic group generated by $[x_1, x_2]N[[F, F], F]$ of order at most 2.

Let $\psi : G_{S'}(2) = F/R \longrightarrow {\rm Gal}(K/k)$ be the natural homomorphism.
Since $[[D_8, D_8], D_8] = 1$ and ramification index of each $\frak{p}_i$ is 2,
$N[[F, F], F] \subset {\rm Ker} (\psi)$. 
Hence $\psi$ induces an isomorphism
\begin{equation*}
F/N[[F, F], F] \stackrel{\sim}{\longrightarrow} {\rm Gal}(K/k).
\end{equation*}
Therefore, $K$ is uniquely determined as a fixed field of $N[[F, F], F]/R$.\\
(2) This follows from Theorem \ref{presentation} (2).
\end{proof}

Suggested by Proposition \ref{millnor residue}, we may define the triple quadratic residue symbol $[\frak{p}_1, \frak{p}_2, \frak{p}_3]$ as follows.

\Def{}{deftriplequadresidue}{
Let $\frak{p}_1 = (\pi_1), \frak{p}_2 = (\pi_2)$ and $\frak{p}_3=(\pi_3)$ be distinct finite primes of $k$ such that ${\rm N}\frak{p}_i = 1$ mod $4$.
Assume that the narrow class group of $k$ is trivial and that
\[\left(\frac{\pi_1}{\pi_2}\right) = \left(\frac{\pi_2}{\pi_1}\right) = 1\ {\rm and}\ \left(\frac{\varepsilon}{\frak{p}_i}\right) = 1\ (1 \leq i \leq 3).\]
Then the triple quadratic residue symbol $[\frak{p}_1, \frak{p}_2, \frak{p}_3]$ is
defined by
\begin{equation*}
[ \frak{p}_1, \frak{p}_2, \frak{p}_3 ] := (-1)^{\mu_2 (123)}.
\end{equation*}
}

\Thm{}{arithmeticmeaning}{
Let the assumptions be as in Definition \ref{deftriplequadresidue}.
Let $K$ be the R\'edei type $D_8$-extension of $k$ for $\{\frak{p}_1, \frak{p}_2\}$ by Theorem \ref{existunique}.
Then we have the following:
\begin{equation*}
[ \frak{p}_1, \frak{p}_2, \frak{p}_3 ] = 1 \Longleftrightarrow \frak{p}_3\  \mbox{is completely decomposed in} \  K/k.
\end{equation*}
}
\begin{proof}
The existence of a R\'edei type $D_8$-extension $K$ of $k$ for $\{\frak{p}_1, \frak{p}_2\}$ follows from Theorem \ref{presentation} (2).
The uniqueness of $K$ follows from Theorem \ref{existunique}. 
The assertion for $[ \frak{p}_1, \frak{p}_2, \frak{p}_3 ]$ follows from Theorem \ref{presentation} (3).
\end{proof}

\Cor{}{}{
Let the assumptions be as in Definition \ref{deftriplequadresidue}.
Then, the following equation holds:
\begin{equation*}
[ \frak{p}_1, \frak{p}_2, \frak{p}_3 ] = [ \frak{p}_2, \frak{p}_1, \frak{p}_3 ].
\end{equation*}
}
\begin{proof}
The construction of R\'edei type $D_8$-extension $K$ of $k$ depends on the set $\{\frak{p}_1, \frak{p}_2\}$ of two primes, and is independent of the ordering of them.
Therefore, by Theorem \ref{arithmeticmeaning}, the values of the two symbols $[ \frak{p}_1, \frak{p}_2, \frak{p}_3 ]$ and $[ \frak{p}_2, \frak{p}_1, \frak{p}_3 ]$ are equal.
\end{proof}

\section{R\'edei type $D_8$-extensions over real quadratic fields}\label{constructredeiextension}
In this section, we give examples of R\'edei type $D_8$-extensions of $k$ for $\{\frak{p_1}, \frak{p}_2\}$ in Theorem \ref{existunique}.
We keep the same notations and assumptions as in the Section \ref{deftriplequadraticresiduesymbol}.
Namely, let $k$ be a real quadratic field whose narrow class number is one.
By Gauss' genus theory (\cite[Remark 4.7, p.172]{onotakashi}), there is a prime number $p = 2$ or $p \equiv 1$ mod $4$ such that $k = \mathbb{Q}(\sqrt{p})$.
We assume that $p \equiv 1$ mod $4$.
Let $\frak{p}_1 = (\pi_1), \frak{p}_2 = (\pi_2), \frak{p}_3 = (\pi_3)$ be distinct finite primes of $k$ such that ${\rm N}\frak{p}_i \equiv 1$ mod 4, $\left(\frac{\varepsilon}{\frak{p}_i}\right) = 1$ ($1 \leq i \leq 3$) and $\left(\frac{\pi_i}{\pi_j}\right) = 1$ ($1 \leq i \neq j \leq 3$).
Let $p_i$ denote the rational prime lying below $\frak{p}_i$.

\Prp{}{inertinertRedei}{
Suppose that both $p_1$ and $p_2$ are inert in $k/\mathbb{Q}$ and $p_1, p_2 \equiv 1$ {\rm mod} $4$.
Let $R$ be R\'edei's extension of $\mathbb{Q}$ for $\{p_1, p_2\}$, which is given by
\begin{equation*}
R = \mathbb{Q}(\sqrt{p_1}, \sqrt{p_2}, \sqrt{x+y\sqrt{p_1}}),
\end{equation*}
where $x, y$ are non-zero integers satisfying
\begin{equation*}
x^2 - p_1 y^2 - p_2 z^2 = 0, \hspace{0.2in} y \equiv 0 \ {\rm mod}\ 2, \hspace{0.2in} x-y \equiv 1 \ {\rm mod}\ 4
\end{equation*}
with some non-zero integer $z$ {\rm (\cite{Redei1939})}.
Then, the R\'edei type $D_8$-extension $K$ of $k$ for $\{\frak{p}_1, \frak{p}_2\}$ is the composite field
\begin{equation*}
K = R \cdot k = k(\sqrt{p_1}, \sqrt{p_2}, \sqrt{x+y\sqrt{p_1}}).
\end{equation*}
}
\begin{proof}
Recall that $R$ is a Galois extension of $\mathbb{Q}$ unramified outside $\{p_1, p_2, \infty\}$ with the ramification index of each $p_i$ being 2 and whose Galois group is isomorphic to the group $D_8$ (\cite{Amano2014}, \cite{Redei1939}).
Since only $\mathbb{Q}(\sqrt{p_1})$, $\mathbb{Q}(\sqrt{p_2})$, $\mathbb{Q}(\sqrt{p_1p_2})$ are the intermediate quadratic fields of $R/\mathbb{Q}$, $R \cap k$ is $\mathbb{Q}$.
By Galois theory, $K/k$ is a Galois extension and Gal($K/k$)$=D_8$.
By Hilbert ramification theory, $K/k$ is unramified outside $\{\frak{p}_1, \frak{p}_2, \infty_1, \infty_2\}$.
Moreover, since $p_1$ and $p_2$ are inert in $k/\mathbb{Q}$, the ramification index of each $\frak{p}_i$ is 2.
Therefore $K$ is the R\'edei type $D_8$-extension of $k$ for $\{\frak{p}_1, \frak{p}_2\}$.
\end{proof}

Next, we consider the case that $p_1$ or $p_2$ is decomposed in $k/\mathbb{Q}$.
To deal with this case, we recall some lemmas regarding ramification.

\Lem{\cite{MR3444843} Kapitel III, (2.5) Theorem, p.209}{relativediscriminantelement}{
Let $L/F$ be an extension of algebraic number fields, $d(L/F)$ be the relative discriminant of $L/F$, and $d(\alpha,L/F)$ be the relative discriminant of $\alpha \in \mathcal{O}_L$ for $L/F$.
Then $d(L/F) \mid (d(\alpha,L/F))$.
}

Using Lemma \ref{relativediscriminantelement}, we obtain the following lemma.

\Lem{}{unramified22}{Let k be a real quadratic field $k = \mathbb{Q}(\sqrt{p})$, $\pi$ be a prime element of $k$ such that $\pi \equiv 1$ {\rm mod} $4\mathcal{O}_k$.
Then the extension $k(\sqrt{\pi})/k$ is ramified at only $(\pi)$.
}
\begin{proof}
We set $\alpha_1 := (1+\sqrt{\pi})/2$.
Since $\alpha_1$ satisfies $\alpha_1^2 - \alpha_1 + (1-\pi)/4 = 0$ and $\pi \equiv 1$ {\rm mod} $4\mathcal{O}_k$, we have $\alpha_1 \in \mathcal{O}_{k(\sqrt{\pi})}$.
Then the relative discriminant of $\alpha_1$ in $k(\sqrt{\pi})/k$ is
\begin{equation*}
d(\alpha_1, k(\sqrt{\pi})/k) = \begin{vmatrix} 1 & \alpha_1 \\ 1 & \overline{\alpha_1} \end{vmatrix}^2 = \pi.
\end{equation*}
By Lemma \ref{relativediscriminantelement}, $k(\sqrt{\pi})/k$ can be ramified at only $(\pi)$.
\end{proof}


\Lem{}{epsilonconjugate}{Let $k$ be a real quadratic field with trivial narrow class group, and let $\pi > 0$ be a prime element of $\mathcal{O}_k$, and let $\varepsilon$ be the fundamental unit of $k$.
We assume $\pi \equiv 1$ {\rm mod} $4\mathcal{O}_k$.
Then, the following are equivalent.\\
{\rm (i)} $\displaystyle \left(\frac{\varepsilon}{\pi}\right) = 1.$\\
{\rm (ii)} $\iota_{\infty_2}(\pi) > 0$, where $\iota_{\infty_2}$ is the conjugate embedding.
}
\begin{proof}
We calculate the Hilbert symbol $\left(\frac{\varepsilon, \pi}{\frak{p}}\right)$ for all primes $\frak{p}$ of $k$ to use product formula.
By definition of the Hilbert symbol, we have
\[
\left(\frac{\varepsilon, \pi}{\frak{p}}\right) = \frac{(\varepsilon, k_\frak{p}(\sqrt{\pi})/k_\frak{p})(\sqrt{\pi})}{\sqrt{\pi}},
\]
where $(\phantom{\varepsilon}, k_\frak{p}(\sqrt{\pi})/k_\frak{p}) : k_\frak{p}^{\times} \longrightarrow {\rm Gal}(k_\frak{p}(\sqrt{\pi})/k_\frak{p})$ is a norm residue symbol.
\underline{Case $\frak{p} \nmid (\pi)\infty$.}
By Lemma \ref{unramified22}, $k_\frak{p}(\sqrt{\pi})/k_\frak{p}$ is an unramified extension.
So, the norm residue symbol factors through ${\rm Gal}(k_\frak{p}^{\rm ur}/k_\frak{p})$ ; 
\[
k_\frak{p}^{\times} \overset{(\phantom{\pi_1}, k_\frak{p}^{\rm ur}/k_\frak{p})} {\xrightarrow{\hspace*{1.5cm}}} {\rm Gal}(k_\frak{p}^{\rm ur}/k_\frak{p}) \overset{{\rm res}} {\xrightarrow{\hspace*{1cm}}} {\rm Gal}(k_\frak{p}(\sqrt{\pi})/k_\frak{p}).
\]
Since the local class field theory, the karnel of $(\phantom{\varepsilon}, k_\frak{p}(\sqrt{\pi})/k_\frak{p})$ is $\mathcal{O}_\frak{p}^{\times}$, and $\varepsilon \in \mathcal{O}_\frak{p}^{\times}$, so $(\varepsilon, k_\frak{p}^{\rm ur}/k_\frak{p}) = {\rm id}_{k_\frak{p}^{\rm ur}}$.
Therefore $(\varepsilon, k_\frak{p}(\sqrt{\pi})/k_\frak{p}) = {\rm id}_{k_\frak{p}(\sqrt{\pi})}$ and hence we have 
$\displaystyle \left(\frac{\varepsilon, \pi}{\frak{p}}\right) = 1$.\\

\underline{Case $\frak{p} = (\pi)$.}
By definition of the quadratic residue symbol, $\displaystyle \left(\frac{\varepsilon, \pi}{\frak{p}}\right) = \left(\frac{\varepsilon}{\pi}\right)$.

\underline{Case $\frak{p} \mid \infty$.}
Since $\frak{p}$ is a real prime, $k_\frak{p} = \mathbb{R}$.
If $\frak{p} = \infty_1$, $(\varepsilon, k_{\frak{p}}(\sqrt{\pi})/k_{\frak{p}}) = (\varepsilon, \mathbb{R}/\mathbb{R}) = {\rm id}_{\mathbb{R}}$ since $\pi > 0$.
Thus, we have $\displaystyle \left(\frac{\varepsilon, \pi}{\frak{p}}\right) = 1$.
If $\frak{p} = \infty_2$, since $(\iota_{\infty_2}(\varepsilon), k_\frak{p}(\sqrt{\iota_{\infty_2}(\pi)})/k_\frak{p}) = (\iota_{\infty_2}(\varepsilon), \mathbb{R}(\sqrt{\iota_{\infty_2}(\pi)})/\mathbb{R})$, then
\[\left(\frac{\varepsilon, \pi}{\frak{p}}\right) = \begin{cases}
1 & {\rm if}\ \iota_{\infty_2}(\pi) > 0 \\
-1 & {\rm if}\ \iota_{\infty_2}(\pi) < 0.
\end{cases}.\]

From the above calculation and Hilbert's product formula, we have
\[
\displaystyle \left(\frac{\varepsilon}{\pi}\right) = 1 \iff \iota_{\infty_2}(\pi) > 0.
\]
\end{proof}

\Thm{}{mainmaintheorem}{Let $k=\mathbb{Q}(\sqrt{p})$ be a real quadratic field with trivial narrow class group and we assume $p \equiv 5\ {\rm mod}\ 8$.
Let $\varepsilon$ be the fundamental unit of $k$, and $\frak{p}_1 = (\pi_1), \frak{p}_2 = (\pi_2)$ are prime ideals of $\mathcal{O}_k$ such that:
\begin{equation*}
\left(\frac{\varepsilon}{\pi_i}\right) = 1,\ \left(\frac{\pi_i}{\pi_j}\right) = 1,\ \pi_i \equiv 1\ {\rm mod}\ 4\mathcal{O}_k,\ \pi_i > 0\ (1 \leq i \neq j \leq 2).
\end{equation*}
Then the followings hold:\\
{\rm (1)} There is a non-zero $\mathcal{O}_k$ solution $(x, y, z)$ of $x^2 - \pi_1 y^2 - \pi_2 z^2 = 0$ such that ${\rm gcd}(x, y, z) = 1$, $y \in 2\mathcal{O}_k$.\\
{\rm (2)} If $x - y \equiv 1\  {\rm mod}\  4\mathcal{O}_k$, the R\'edei type $D_8$-extension $K$ of $k$ for $\{\frak{p}_1, \frak{p}_2\}$ is given by
\begin{equation*}
K = k(\sqrt{\pi_1}, \sqrt{\pi_2}, \sqrt{x+y\sqrt{\pi_1}}).
\end{equation*}
}
\begin{proof}
(1) We use Hasse--Minkowski theorem to show the equation in (1) has a nontrivial $\mathcal{O}_k$ solution.
It suffices to show that the Hilbert symbol $\displaystyle \left(\frac{\pi_1, \pi_2}{\frak{p}}\right)$ is $1$ for all primes $\frak{p}$ of $k$.\\
\underline{Case $\frak{p} \nmid \frak{p}_1\frak{p}_2\infty$.}
Since $\pi_2 \in \mathcal{O}_\frak{p}^{\times}$ and Lemma \ref{unramified22}, $k_\frak{p}(\sqrt{\pi_2})/k_\frak{p}$ is an unramified extension.
Furthermore, like the proof of Lemma \ref{epsilonconjugate}, it can be seen that $\pi_1 \in \mathcal{O}_\frak{p}^{\times}$, $(\pi_1, k_\frak{p}(\sqrt{\pi_2})/k_\frak{p}) = {\rm id}_{k_\frak{p}(\sqrt{\pi_2})}$, and hence we have $\displaystyle \left(\frac{\pi_1, \pi_2}{\frak{p}}\right) = 1$.\\
\underline{Case $\frak{p} \mid \frak{p}_1\frak{p}_2$.}
Since the quadratic residue symbols $\displaystyle \left(\frac{\pi_2}{\pi_1}\right) = \left(\frac{\pi_1}{\pi_2}\right) = 1$, we have $\displaystyle \left(\frac{\pi_1, \pi_2}{\frak{p}}\right) = 1$.\\
\underline{Case $\frak{p} \mid \infty$.}
Since $\frak{p}$ is a real prime, $k_\frak{p} = \mathbb{R}$.
If $\frak{p} = \infty_1$, $(\pi_1, k_\frak{p}(\sqrt{\pi_2})/k_\frak{p}) = (\pi_1, \mathbb{R}(\sqrt{\pi_2})/\mathbb{R}) = (\pi_1, \mathbb{R}/\mathbb{R}) =1$ and hence we have $\displaystyle \left(\frac{\pi_1, \pi_2}{\frak{p}}\right) = 1$.\\
If $\frak{p} = \infty_2$, $\iota_{\infty_2}(\pi_2) > 0$ by Lemma \ref{epsilonconjugate}.
Therefore $(\pi_1, k_\frak{p}(\sqrt{\iota_{\infty_2}(\pi_2)})/k_\frak{p}) = (\pi_1, \mathbb{R}(\sqrt{\iota_{\infty_2}(\pi_2)})/\mathbb{R}) = (\pi_1, \mathbb{R}/\mathbb{R}) =1$ and hence we have $\displaystyle \left(\frac{\pi_1, \pi_2}{\frak{p}}\right) = 1$.\\
From the above, the Hilbert symbol $\displaystyle \left(\frac{\pi_1, \pi_2}{\frak{p}}\right)$ is $1$ for all primes $\frak{p}$ of $k$.
Then there is a non-zero $\mathcal{O}_k$ solution of $x^2 - \pi_1 y^2 - \pi_2 z^2 = 0$ by Hasse--Minkowski theorem.
Since the class number of $k$ is one, we can assume ${\rm gcd}(x, y, z) = 1$.
Since $x^2 - \pi_1 y^2 - \pi_2 z^2 = 0$, we have $x^2 \equiv y^2 + z^2$ {\rm mod} $4\mathcal{O}_k$.
If $x^2 \equiv 0$ {\rm mod} $4\mathcal{O}_k$, then it causes a contradiction to ${\rm gcd}(x, y, z) = 1$.
Therefore we can set $y^2 \equiv 0$ {\rm mod} $4\mathcal{O}_k$.
This means that $y \in 2\mathcal{O}_k$.\\
(2) First, we easily see that $K = k(\sqrt{\pi_1}, \sqrt{\pi_2}, \sqrt{x+y\sqrt{\pi_1}})$ is a $D_8$-extension of $k$, where $\alpha_1 := x + y \sqrt{\pi_1}$.
In fact, to be precise, the intermediate fields of $K$ over $k$ are given in the following, where $k_i = k(\sqrt{\pi_i})$, and $\alpha_2 = (\sqrt{\alpha_1} + \sqrt{\bar{\alpha}_1})^2$.
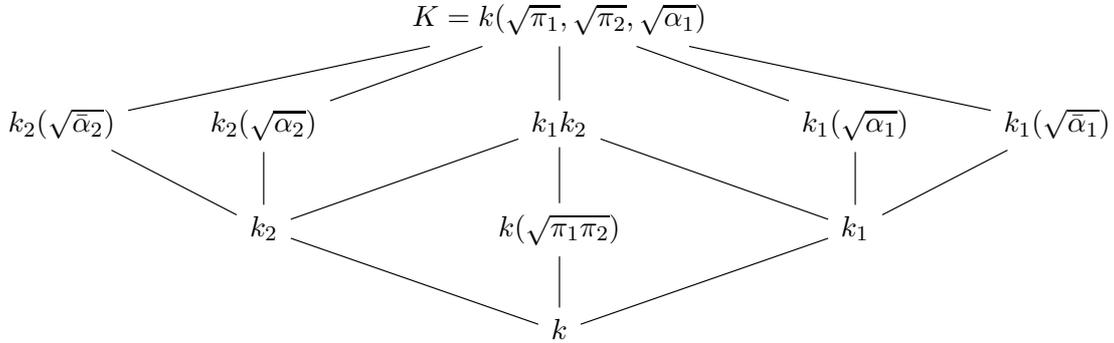
\begin{equation*}
\begin{tikzcd}
   & & K = k(\sqrt{\pi_1},\sqrt{\pi_2}, \sqrt{\alpha_1})\\
 k_2(\sqrt{\bar{\alpha}_2}) \arrow[rru, -] & k_2(\sqrt{\alpha_2}) \arrow[ru, dash] & k_1k_2 \arrow[u, dash] & k_1(\sqrt{\alpha_1}) \arrow[lu, dash] & k_1(\sqrt{\bar{\alpha}_1}) \arrow[llu, dash] \\
 & k_2 \arrow[lu, dash] \arrow[u, dash] \arrow[ru, dash] & k(\sqrt{\pi_1\pi_2}) \arrow[u, dash] & k_1 \arrow[lu, dash] \arrow[u, dash] \arrow[ru, dash] \\
 && k \arrow[lu, dash] \arrow[u, dash] \arrow[ru, dash]
 \end{tikzcd}
\end{equation*}
Next, let us see the ramifications in $K/k$.
We prove that $k_1(\sqrt{\alpha_1}) / k_1$ is ramified at only primes dividing $\pi_2$.
We first claim that $\theta := \frac{1+\sqrt{\alpha_1}}{2} \in \mathcal{O}_{k_1(\sqrt{\alpha_1})}$.
We set $\beta := \frac{1-\alpha_1}{4}$.
It can be seen that $\beta$ satisfies $\beta^2 - (\frac{1-(x-y)}{2} - \frac{y}{2}) \beta + \frac{1}{16}((1-x)^2 - y^2 + (1-\pi_1) y^2) = 0$.
Since $x-y \equiv 1$ mod $4\mathcal{O}_k$ and $y \in 2\mathcal{O}_k$, $1-x+y \in 4\mathcal{O}_k$ and $1-x-y = 1-x+y-2y \in 4\mathcal{O}_k$.
Then $(1-x)^2 - y^2 = (1-x+y)(1-x-y) \in 16\mathcal{O}_k$.
Since $\pi_1 \equiv 1$ mod $4\mathcal{O}_k$ and $y \in 2\mathcal{O}_k$, $(1-\pi_1) y^2 \in 16\mathcal{O}_k$.
Therefore, $\beta \in \mathcal{O}_{k_1}$.
Since $\theta$ satisfies $\theta^2 - \theta + \beta = 0$, $\theta \in \mathcal{O}_{k_1(\sqrt{\alpha_1})}$.
Since the relative discriminant of $\theta$ in $k_1(\sqrt{\alpha_1}) / k_1$ is
\begin{equation*}
d(\theta, k_1(\sqrt{\alpha_1}) / k_1) = \begin{vmatrix} 1 & \theta \\ 1 & \overline{\theta} \end{vmatrix}^2 = \alpha_1.
\end{equation*}
Since ${\rm N}_{k_1(\sqrt{\alpha_1}) / k_1}(\alpha_1) = \pi_2$ and Lemma \ref{relativediscriminantelement}, $k_1(\sqrt{\alpha_1}) / k_1$ is ramified at only primes dividing $\pi_2$.
\end{proof}

Here is a numerical example of Theorem \ref{mainmaintheorem}.

\Ex{}{exampleborromean1}{
Let $k = \mathbb{Q}(\sqrt{5})$, $\pi_1 = 33 + 8\sqrt{5}\ ({\rm N}_{k/\mathbb{Q}}(\pi_1) = 769)$, $\pi_2 = 17$, $\pi_3 = \frac{23+5\, \sqrt{5}}{2}\ ({\rm N}_{k/\mathbb{Q}}(\pi_3) = 101)$.
We set $\frak{p}_i = (\pi_i)\ (1 \leq i \leq 3)$.
Then we have $$\displaystyle \left(\frac{\varepsilon}{\pi_i}\right) = 1,\ \left(\frac{\pi_i}{\pi_j}\right) = 1,\ \pi_i \equiv 1\ {\rm mod}\ 4\mathcal{O}_k,\ \pi_i > 0\ (1 \leq i \neq j \leq 2),$$
\begin{equation*}
\left(\frac{\varepsilon}{\pi_3}\right) = 1,\ \left(\frac{\pi_i}{\pi_j}\right) = 1,\ {\rm N}\frak{p}_3 \equiv 1\ {\rm mod}\ 4\ (1 \leq i \neq j \leq 3),
\end{equation*}
where $\varepsilon = (1+\sqrt{5})/2$.
It is easy to see that
\begin{center}
$(-23-14\sqrt{5})^2 - \pi_1 2^2 - \pi_2 (2+3\sqrt{5})^2 = 0$ and $(-23-14\sqrt{5}) - 2 \equiv 1\ {\rm mod}\ 4\mathcal{O}_k$.\\
\end{center}
Therefore, $K = k(\sqrt{\pi_1},\sqrt{\pi_2}, \sqrt{\alpha_1})$ is the R\'edei type $D_8$-extension of $k$ for $\{\frak{p}_1, \frak{p}_2\}$, where $\alpha_1 = (-23-14\sqrt{5}) + 2 \sqrt{\pi_1}$.
Since $\frak{p}_3$ is not completely decomposed in $K$, we have
\begin{equation*}
[\frak{p}_1, \frak{p}_2, \frak{p}_3] = -1.
\end{equation*}
}

\Rem{}{}{From the viewpoint of arithmetic topology, mod 2 arithmetic Milnor invariants may be regarded as an arithmetic analogue of (higher) linking numbers of a link.
So, the triple of primes $(\frak{p}_1 = (33+8\sqrt{5}), \frak{p}_2 = (17), \frak{p}_3 = (\frac{23+5\sqrt{5}}{2}))$ in Example \ref{exampleborromean1} may be called the {\it Borromean primes}.
}

Next, we give other numerical examples, which are not contained in the cases of Proposition \ref{inertinertRedei} and Theorem \ref{mainmaintheorem}.

\Ex{}{decomposedinertRedei}{The case that $p_1$ is decomposed and $p_2$ is inert in $k/\mathbb{Q}$.
We give an example for $k = \mathbb{Q}(\sqrt{5})$, $p_1 = 29$ and $p_2 = 13$. 
We take $\pi_1 = (11+\sqrt{5})/2$ and $\pi_2 = 13$.
Then $\left(\frac{\pi_1}{\pi_2}\right) = \left(\frac{\pi_2}{\pi_1}\right) = 1$
and $\left(\frac{\varepsilon}{\frak{p}_1}\right) = \left(\frac{\varepsilon}{\frak{p}_2}\right) = 1$, where $\varepsilon = (1+\sqrt{5})/2$.
In this case, $K = k(\sqrt{\pi_1},\sqrt{\pi_2}, \sqrt{\alpha_1})$ is the R\'edei type $D_8$-extension of $k$ for $\{\frak{p}_1, \frak{p}_2\}$, where $\alpha_1 = (-1-9\sqrt{5} + 6\sqrt{\pi_1})/4$.

\begin{proof}
Let $k_1 := k(\sqrt{\pi_1})$.
We prove that $k_1/k$ can be ramified at only over $\frak{p}_1$ and $k_1(\sqrt{\alpha_1})/k_1$ is ramified at only primes dividing $\pi_2$.
First we set $\lambda_1 := (1+\sqrt{5}+2\sqrt{\pi_1})/4$.
Since $\lambda_1$ satisfies $\lambda_1^2 - \frac{1+\sqrt{5}}{2} \lambda_1 - 1 = 0$, $\lambda_1 \in \mathcal{O}_{k_1}$.
Since the relative discriminant of $\lambda_1$ in $k_1/k$ is
\begin{equation*}
d(\lambda_1, k_1/k) = \begin{vmatrix} 1 & \lambda_1 \\ 1 & \overline{\lambda_1} \end{vmatrix}^2 = \pi_1.
\end{equation*}
By Lemma \ref{relativediscriminantelement}, $k_1/k$ can be ramified at only over $\frak{p}_1$.
Next we set $\lambda_2 := ((\sqrt{5}-1) - (\sqrt{5}+1)\sqrt{\pi_1} + 4\sqrt{\alpha_1})/8$.
Since $\lambda_2$ satisfies $\lambda_2^8 + \lambda_2^7 -2\lambda_2^6 +3\lambda_2^5 +11\lambda_2^4 -\lambda_2^3 -18\lambda_2^2 -13\lambda_2 -1 = 0$, $\lambda_2 \in \mathcal{O}_{k_1(\sqrt{\alpha_1})}$.
Since the relative discriminant of $\lambda_2$ in $k_1(\sqrt{\alpha_1})/k_1$ is
\begin{equation*}
d(\lambda_2, k_1(\sqrt{\alpha_1})/k_1) = \begin{vmatrix} 1 & \lambda_2 \\ 1 & \overline{\lambda_2} \end{vmatrix}^2 = \alpha_1.
\end{equation*}
Since ${\rm N}_{k_1/k}(\alpha_1) = \pi_2$ and Lemma \ref{relativediscriminantelement}, $k_1(\sqrt{\alpha_1})/k_1$ is ramified at only primes dividing $\pi_2$.
Therefore the extension $K/k$ is the R\'edei type $D_8$-extension of $k$ for $\{\frak{p}_1, \frak{p}_2\}$.
\end{proof}}

\Ex{}{examplen}{The case that $p_1$ and $p_2$ are decomposed in $k/\mathbb{Q}$.
We give an example $k = \mathbb{Q}(\sqrt{5})$, $p_1=29$ and $p_2=89$.
We take $\pi_1 = (11+\sqrt{5})/2$ and $\pi_2 = (19+\sqrt{5})/2$.
Then $\left(\frac{\pi_1}{\pi_2}\right) = \left(\frac{\pi_2}{\pi_1}\right) = 1$ and $\left(\frac{\varepsilon}{\frak{p}_1}\right) = \left(\frac{\varepsilon}{\frak{p}_2}\right) = 1$, where $\varepsilon = (1+\sqrt{5})/2$.
In this case, $K = k(\sqrt{\pi_1},\sqrt{\pi_2}, \sqrt{\alpha_1})$ is the R\'edei type $D_8$-extension of $k$ for $\{\frak{p}_1, \frak{p}_2\}$, where $\alpha_1 = (6\sqrt{5} + (1-\sqrt{5}) \sqrt{\pi_1})/4$.
\begin{proof}
Let $k_1 := k(\sqrt{\pi_1})$ and $k_2 := k(\sqrt{\pi_2})$.
We prove that $k_2/k$ is ramified at only $\frak{p}_2$ and $k_1(\sqrt{\alpha_1})/k_1$ is ramified at only primes dividing $\pi_2$.
First we set $\theta_1 := (1+\sqrt{5}+2\sqrt{\pi_2})/4$.
Since $\theta_1$ satisfies $\theta_1^2 - \frac{1+\sqrt{5}}{2} \theta_1 - 2 = 0$, $\theta_1 \in \mathcal{O}_{k_2}$.
Since the relative discriminant of $\theta_1$ in $k_2/k$ is
\begin{equation*}
d(\theta_1, k_2/k) = \begin{vmatrix} 1 & \theta_1 \\ 1 & \overline{\theta_1} \end{vmatrix}^2 = \pi_2.
\end{equation*}
By Lemma \ref{relativediscriminantelement}, $k_2/k$ is ramified at only over $\frak{p}_2$.
Next we set $\theta_2 := ((2\sqrt{5}+4) - (\sqrt{5}+3)\sqrt{\pi_1} + 4\sqrt{\alpha_1})/8$.
Since $\theta_2$ satisfies $\theta_2^8 -4 \theta_2^7 +3\theta_2^5 -2\theta_2^4 +12\theta_2^3 +2\theta_2^2 +5\theta_2 -1 = 0$, $\theta_2 \in \mathcal{O}_{k_1(\sqrt{\alpha_1})}$.
Since the relative discriminant of $\theta_2$ in $k_1(\sqrt{\alpha_1})/k_1$ is
\begin{equation*}
d(\theta_2, k_1(\sqrt{\alpha_1})/k_1) = \begin{vmatrix} 1 & \theta_2 \\ 1 & \overline{\theta_2} \end{vmatrix}^2 = \alpha_1.
\end{equation*}
Since ${\rm N}_{k_1/k}(\alpha_1) = \pi_2$ and Lemma \ref{relativediscriminantelement}, $k_1(\sqrt{\alpha_1})/k_1$ is ramified at only primes dividing $\pi_2$.
Therefore the extension $K/k$ is the R\'edei type $D_8$-extension of $k$ for $\{\frak{p}_1, \frak{p}_2\}$.
\end{proof}}

\section{Massey products in Galois cohomology}\label{masseymassey}
In this section, we interpret our mod 2 arithmetic Milnor invariants and triple quadratic residue symbols in terms of the Massey product in Galois cohomology.
Our theorem is seen as a generalization of the known relation between the cup product and the quadratic residue symbol to the triple case, and also a generalization of the previous result by Morishita (\cite{Morishita2004}) in the case of the rational number field to the real quadratic fields.
It may be regarded as a mod 2 arithmetic analogue of the corresponding topological result due to Turaev (\cite{Turaev1979}).\\

Let $G$ be a pro-$2$ group and let $A$ be a commutative ring with identity on which $G$ acts trivally.
Let $C^j(G,A)$ be the $A$-module of inhomogeneous $j$-cochains $(j \geq 0)$ of $G$ with coefficients in $A$ and we consider the differential graded algebra $(C^{\bullet}(G,A), d)$, 
where the product structure on $C^{\bullet}(G,A) = \bigoplus_{j\geq 0} C^j(G,A)$ is given by the cup product and the differential $d$ is the coboundary operator.
Then we have the cohomology $H^*(G,A) = H^*(C^{\bullet}(G,A))$ of the pro-$2$ group $G$ with coefficients in $A$.
In the following, we consider Massey products in $H^1(G,A)$.
For the sign convention, we follow \cite{Dwyer1975}.
Let $\chi_1,\dots,\chi_n \in H^1(G,A)$ $(n\geq 2)$.
An $n$-th {\it Massey product}  $\langle \chi_1,\dots,\chi_n \rangle$ is said to be {\it defined} if there is an array
\begin{equation*}
\Omega = \{\omega_{ij} \in C^1(G,A) \; | \; 1 \leq i < j \leq n+1, (i,j) \neq (1,n+1) \}
\end{equation*}
such that
\begin{equation*}
\left\{ 
\begin{array}{l}
[\omega_{i,i+1}] = \chi_i \;\; (1\leq i \leq n),\\
 \displaystyle{d\omega_{ij} = \sum_{a=i+1}^{j-1} \omega_{ia}\cup \omega_{aj}} \;\; (j \neq i+1).
\end{array}
\right.
\end{equation*}
Such an array $\Omega$ is called a {\it defining system} for $\langle \chi_1,\dots, \chi_n \rangle$. The value of  $\langle \chi_1,\dots, \chi_n \rangle$ relative to $\Omega$, denoted by $\langle \chi_1,\dots, \chi_n \rangle_{\Omega}$, is defined by the cohomology class represented by the $2$-cocycle
\begin{equation*}
\sum_{a=2}^n \omega_{1a} \cup \omega_{a,n+1}.
\end{equation*}
We define the Massey product $\langle \chi_1,\dots,\chi_n \rangle $ to be the subset of $H^2(G,A)$ consisting of elements $\langle \chi_1,\dots, \chi_n \rangle_{\Omega}$ for some defining system $\Omega$.
By convention, $\langle \chi \rangle = 0$.
We recall the following basic fact (cf. \cite{Kraines1966}).
\Lem{}{Masseyproductlem}{We have $\langle \chi_1, \chi_2 \rangle = \chi_1 \cup \chi_2$.
For $n \geq 3$, $\langle \chi_1,\dots, \chi_n \rangle$ is defined and consists of a single element if $\langle \chi_{j_1},\dots, \chi_{j_a} \rangle = 0$ for all proper subsets $\{j_1,\dots, j_a \}$ $(a \geq 2)$ of $\{1,\dots ,n\}$.
(In this case, we denote the single element by $\langle \chi_{1},\dots, \chi_{n} \rangle$.)}

Next, we recall a relation between Massey products and Magnus coefficients.
Suppose that $G$ is a finitely generated pro-$2$ group with a minimal presentation 
\begin{equation*}
1 \longrightarrow N \longrightarrow F \stackrel{\psi}{\longrightarrow} G \longrightarrow 1,
\end{equation*}
where $F$ is a free pro-$2$ group on generators $x_1, \dots, x_n$ with $n = {\rm dim}_{\mathbb{F}_2} H^1(G, \mathbb{F}_2)$.
We set $\tau_i := \psi(x_i) \;(1\leq i \leq n)$.
We assume that $\psi$ induces the isomorphism $F/\Phi(F) \simeq G/\Phi(G)$ so that $\psi$ induces the isomorphism $\psi^* : H^1(G, \mathbb{F}_2) \simeq H^1(F, \mathbb{F}_2)$.
We let 
\begin{equation*}
{\rm tra} : H^1(N, \mathbb{F}_2)^{G} \rightarrow H^2(G, \mathbb{F}_2)
\end{equation*}
 be the transgression defined as follows.
 For $a \in H^1(N, \mathbb{F}_2)^{G}$, choose a 1-cochain $b \in C^1(F, \mathbb{F}_2)$ such that $b|_{N} = a$.
Since the value $db(f_1,f_2)$, $f_i \in F$, depends only on the cosets $f_i$ mod $N$, $db$ defines a 2-cocyle $c$ of $G$.
Then ${\rm tra}(a)$ is defined by the class of $c$.
By the Hochschild-Serre spectral sequence, ${\rm tra}$ is an isomorphism and so we have the dual isomorphism, called the Hopf isomorphism,
\begin{equation*}
{\rm tra}^{\vee} : H_2(G, \mathbb{F}_2) \stackrel{\sim}{\rightarrow} H_1(N, \mathbb{F}_2)_{G} = N/N^2[N, F].
\end{equation*}
Then we have the following proposition (cf. \cite[Lemma 1.5]{Stein1990}).
The proof goes in the same manner as in \cite[Theorem 2.2.2]{Morishita2004}.
\Prp{}{Masseyproductprop}{Notations being as above, let $\chi_1,\dots,\chi_n \in H^1(G, \mathbb{F}_2)$ $(n \geq 2)$.
Let $f \in N$ and set $\delta := ({\rm tra}^{\vee})^{-1}(f \; {\rm mod}\; N^2[N,F])$.
Assume that all Massey products up to length $n-1$ are trivial.
Then $N \subset F^{(n)}$ and we have
\begin{equation*}
\langle \chi_1,\dots, \chi_n \rangle(\delta) = \sum_{{\scriptstyle I = (i_1\cdots i_n)} \atop {\scriptstyle 1 \leq i_1,\dots , i_n \leq N}} \chi_1(x_{i_1})\cdots \chi_n(x_{i_n})\mu_l(I; f).
\end{equation*}}

Let us be back in our arithmetic situation and keep the same notations as in Section \ref{deftriplequadraticresiduesymbol}.
So let $k$ be the real quadratic field and let $S := \{ \frak{p}_1, \frak{p}_2, \frak{p}_3\}$ be a set of all distinct primes where ${\rm N}\frak{p}_i \equiv 1$ mod $4$ and $\left(\frac{\varepsilon}{\frak{p}_i}\right) = 1$.
Let $\overline{S} := S \cup S_{\mathbb{R}}$ where $S_{\mathbb{R}} = \{\infty_1, \infty_2\}$ be the set of real primes of $k$.
By Theorem \ref{maintheorem1}, we have the following minimal presentation of the Galois group $G_S(2)$ of maximal pro-$2$ extension over $k$ unramified outside $\overline{S}$
\begin{align*}
G_S(2) &= \left \langle x_1, x_2 , x_3 \, \middle | \begin{array}{c} 
x_1^{\mathrm{N}\frak{p}_{1}-1} [ x_1, y_1] = x_2^{\mathrm{N}\frak{p}_{2}-1} [ x_2, y_2] = x_3^{\mathrm{N}\frak{p}_3-1} [ x_3, y_3] =1\\
\\
x_{\infty_2}^2 = 1
\end{array} \right \rangle \\
&= F_S/N_S.  \end{align*}
Here $x_i$ denotes the word representing a monodromy $\tau_i$ over $\frak{p}_i$ in $k_S(2)/k$ $(1\leq i \leq 3)$ and $F_S$ denotes the free pro-$2$ group on $x_1, x_2, x_3$.
The pro-$2$ word $y_i$ represents a Frobenius automorphism $\sigma_i$ over $\frak{p}_i$ in $k_S(2)/k$ and  $N_S$ denotes the closed subgroup of $F_S$ generated normally by $x_{i}^{{\rm N}\frak{p}_{i} -1}[x_{i},y_{i}]$ for $1\leq i \leq 3$ and $x_{\infty_2}^2$.
We set $\delta_i := ({\rm tra}^{\vee})^{-1}(x_i^{{\rm N}\frak{p}_i -1}[x_i,y_i])$, where ${\rm tra}^{\vee} : H_2(G_S(2),\mathbb{F}_2) \stackrel{\sim}{\rightarrow} N_S/N_S^2[N_S,F_S]$ is the Hopf isomorphism.
Let $\chi_1, \chi_2, \chi_3 \in H^1(G_S(2),\mathbb{F}_2)$ be the Kronecker dual to the monodromies $\tau_1, \tau_2, \tau_3$, namely, $\chi_i(\tau_j) = \delta_{i,j}$.
\Thm{}{}{Notations being as above, we have, for $1\leq i \neq j \leq 3$,
\begin{equation*}
\left(\frac{\pi_i}{\pi_j} \right) = (-1)^{ \langle \chi_i, \chi_j \rangle(\delta_j)},\;\; \left(\frac{\pi_j}{\pi_i}\right) = (-1)^{\langle \chi_i, \chi_j \rangle(\delta_i)}.
\end{equation*}}
\begin{proof}
By Proposition \ref{Masseyproductprop}, we have
$$\langle \chi_i, \chi_j \rangle(\delta_a) = \mu_2((ij); x_a^{{\rm N}\frak{p}_a-1}[x_a,y_a]).$$
Here, noting that the Magnus coefficients are in $\mathbb{F}_2$, we have
\begin{gather*}
M(x_a^{{\rm N}\frak{p}_a-1}[x_a,y_a]) = M(x_a^{{\rm N}\frak{p}_a-1})M(x_a)M(y_a)M(x_a)^{-1}M(y_a)^{-1},\\
M(x_a) = 1 + X_a, M(x_a^{{\rm N}\frak{p}_a-1}) = 1 + ({\rm terms\  of\  degree  \geq 4}),\\
M(x_a)^{-1} = 1 + X_a + X_a^2 + ({\rm terms\  of\  degree  \geq 3}),\\
M(y_a) = 1 + \sum_i \mu_2(ia)X_i + \sum_{i, j} \mu_2(ija)X_iX_j + ({\rm terms\  of\  degree  \geq 3}),\\
M(y_a)^{-1} = 1 + \sum_i \mu_2(ia)X_i + \sum_{i, j} \{\mu_2(ija) + \mu_2(ia)\mu_2(ja)\}X_iX_j + ({\rm terms\  of\  degree  \geq 3}),
\end{gather*}
The straightforward calculation yields
\begin{equation*}
M(x_a^{{\rm N}\frak{p}_a-1}[x_a,y_a]) = 1+\sum_i \mu_2(ia)X_aX_i + \sum_i \mu_2(ia)X_iX_a + ({\rm terms\  of\  degree  \geq 3}).
\end{equation*}
Hence we obtain
\begin{equation*}
\langle \chi_i, \chi_j \rangle(\delta_a) = \left\{
  \begin{array}{ll}
  \mu_2(ij) & a = j,\\
  \mu_2(ji) & a = i,\\
  0 & a \neq i, j.
  \end{array} \right.
\end{equation*}
The assertion follows from Theorem \ref{millnor residue}.
\end{proof}
We assume that 
\begin{equation*}
\langle \chi_i, \chi_j \rangle  = 0 \;\; (1\leq i \neq j \leq 3),
\end{equation*}
which is equivalent to the condition
\begin{equation*}
\left( \frac{\pi_i}{\pi_j} \right) = 1 \;\; (1\leq i \neq j \leq 3)
\end{equation*}
by Proposition \ref{invarianttheorem}, and the mod $2$ Milnor invariants $\mu_2(abc)$ ($\{a,b,c\} = \{1,2,3\}$) are well defined.
By the definition of Massey products and Lemma \ref{Masseyproductlem}, there are 1-cochains $\omega_{13}, \omega_{24} \in C^1(G_S(2), \mathbb{F}_2)$ such that
\begin{equation*}
\langle \chi_1, \chi_2 \rangle = d\omega_{13}, \; \langle \chi_2, \chi_3 \rangle = d\omega_{24},
\end{equation*}
and we have the triple Massey product $\langle \chi_1, \chi_2, \chi_3 \rangle$ defined by
\begin{equation*}
\langle \chi_1, \chi_2, \chi_3 \rangle = [\chi_1 \cup \omega_{24} + \omega_{13} \cup \chi_3].
\end{equation*}

\Thm{}{}{Assume that
\begin{equation*}
\left( \frac{\pi_i}{\pi_j} \right) = 1 \;\; (1\leq i \neq j \leq 3).
\end{equation*}
Then we have
\begin{equation*}[\frak{p}_1, \frak{p}_2, \frak{p}_3] = (-1)^{\langle \chi_1, \chi_2, \chi_3 \rangle(\delta_3)}, [\frak{p}_2, \frak{p}_3, \frak{p}_1] = (-1)^{\langle \chi_1, \chi_2, \chi_3 \rangle(\delta_1)}.
\end{equation*}
In particular, $\mu(123)$ and $[\frak{p}_1, \frak{p}_2, \frak{p}_3]$ are independent of choices of 1-cochains $\omega_{13}$ and $\omega_{24}$.}
\begin{proof}
Noting by the assumption that all $\mu_2(ij) = 0$, we have
\begin{gather*}
M(x_a) = 1 + X_a, 
M(x_a^{{\rm N}\frak{p}_a-1}) = 1 + ({\rm terms\  of\  degree  \geq 4}),\\
M(x_a)^{-1} = 1 + X_a + X_a^2 + X_a^3 + ({\rm terms\  of\  degree  \geq 4}),\\
M(y_a) = \displaystyle{1 + \sum_{i, j} \mu_2(ija)X_iX_j + \sum_{i, j, k} \mu_2(ijka)X_iX_jX_k} + ({\rm terms\  of\  degree  \geq 4}), \\
M(y_a)^{-1} = \displaystyle{1+ \sum_{i, j} \mu_2(ija) X_iX_j + \sum_{i, j, k} \mu_2(ijka)X_iX_jX_k} + ({\rm terms\  of\  degree  \geq 4}).
\end{gather*}
The straightforward calculation yields
\begin{equation*}
M(x_a^{{\rm N}\frak{p}_a-1}[x_a,y_a]) = 1+\sum_{i, j} \mu_2(ija)X_aX_iX_j + \sum_{i, j} \mu_2(ija)X_iX_jX_a + ({\rm terms\  of\  degree  \geq 4}).
\end{equation*}
Hence we obtain
\begin{equation*}
\langle \chi_1, \chi_2, \chi_3 \rangle(\delta_a) = \left\{
  \begin{array}{ll}
  \mu_2(123) & a = 3,\\
  \mu_2(231) & a = 1,\\
  0 & a \neq 1, 3.
  \end{array} \right.
\end{equation*}
The assertion follows from Defnition \ref{deftriplequadresidue}.
\end{proof}

\subsection*{Acknowledgements}
I would like to express my sincere gratitude to my supervisor, Professor Masanori Morishita, for suggesting the problem studied in this paper and for his guidance and cooperation in carrying out this research.
I am grateful to Professors Yasushi Mizusawa, Dohyeong Kim and Doctor Sosuke Sasaki for useful communications. 
I would like to thank Yuki Ishida and Dingchuan Zheng for useful discussions and suggestions.


\end{document}